\documentclass[12pt]{article}
\usepackage{amsfonts, amsmath}
\newtheorem{theo}{Theorem}[section]
\newtheorem{pro}[theo]{Proposition}
\newtheorem{lem}[theo]{Lemma}
\newtheorem{cor}[theo]{Corollary}

\newcommand{\ccc}{{\cal C}}

\def\norm#1{\left \|\, #1 \, \right \|}

\def \nat{ { {\rm I}\!{\rm N}} }
 
 \newcommand{\ep}{\varepsilon}

\title{Riesz transform on manifolds  and  heat kernel regularity}

\author{Pascal Auscher\thanks{Research partially supported  by the European Commission (IHP
Network ``Harmonic
Analysis and Related Problems'' 2002-2006, Contract HPRN-CT-2001-00273-HARP).} \\Universit\'e de Paris-Sud
\\pascal.auscher@math.u-psud.fr
\and
Thierry Coulhon\thanks{
Research partially supported by  Macquarie University and by the European Commission (IHP
Network ``Harmonic
Analysis and Related Problems'' 2002-2006, Contract HPRN-CT-2001-00273-HARP).} \\
Universit\'e de Cergy-Pontoise\\
Thierry.Coulhon@math.u-cergy.fr 
\and  
Xuan Thinh Duong\thanks{
Research partially supported by the ARC.}\\
Macquarie University, Sydney\\
duong@ics.mq.edu.au
\and  
Steve Hofmann\\
University of Missouri, Columbia\\
hofmann@math.missouri.edu
}
\date{revised; October  25, 2004}

\begin{document}

\maketitle

\abstract{ 
On consid\`ere la classe des vari\'et\'es riemanniennes compl\`etes non
compactes dont le noyau de la chaleur
satisfait une estimation sup\'erieure et inf\'erieure gaussienne.
On montre que la transform\'ee de Riesz y est born\'ee sur $L^p$,
pour un intervalle ouvert de $p$ au-dessus de $2$, si et seulement si
le gradient du noyau de la chaleur satisfait une certaine estimation
$L^p$
pour  le m\^eme intervalle d'exposants $p$.

One considers the class of complete non-compact Riemannian manifolds
whose heat kernel
satisfies Gaussian estimates from above and below. One shows that the
Riesz transform is $L^p$ bounded on such a manifold,
for $p$ ranging in an open interval above $2$, if and only if the
gradient of the heat kernel satisfies a certain $L^p$ estimate
in the same interval of $p$'s.

{\bf MSC numbers 2000:} 58J35, 42B20
}

\tableofcontents

\section{Introduction}

The aim of this article is to give a necessary and sufficient condition for the 
two natural definitions of homogeneous first order  $L^p$ Sobolev   spaces
to coincide on a large class of  Riemannian manifolds, for $p$ in an   interval $(q_0,p_0)$, where $2<p_0\le\infty$
and $q_0$ is the conjugate exponent to  $p_0$. On closed manifolds, these definitions are well-known to coincide for all
$1<p<\infty$.  For
non-compact manifolds, and again $p_0=\infty$,  a sufficient condition    has been asked for by Robert Strichartz in 1983
(\cite{St}) and many partial answers have been given since. We shall review them in Section
\ref{history} below. 
The condition we propose is in terms
of   regularity  of the heat kernel, more precisely in terms of integral estimates of its gradient. We are able to treat
manifolds with the doubling property  together with natural heat kernel bounds,  as well as the ones with locally bounded
geometry where the bottom of the spectrum of the Laplacian is positive.

\subsection{Background}

Let $M$ be a complete  non-compact connected Riemannian manifold,  $\mu$ the  Riemannian measure, $\nabla$  the
Riemannian gradient. 
Denote by $|.|$ the length in the tangent space,
and by $\|.\|_p$ the norm in $L^p(M,\mu)$, $1\le p\le \infty$. One defines $\Delta$, the
Laplace-Beltrami operator, as a self-adjoint positive operator on $L^2(M,\mu)$ by the formal
integration by parts 
\begin{equation*}
(\Delta f,f)= \||\nabla f|\|_2^2
\end{equation*} 
for all 
$f\in\ccc^\infty_0(M)$, and its positive self-adjoint square root $\Delta^{1/2}$ by
\begin{equation*}
(\Delta f,f)=\|\Delta^{1/2}f\|_2^2.
\end{equation*}  As a consequence,
\begin{equation}
\||\nabla f|\|_2^2=\|\Delta^{1/2}f\|_2^2. \label{triv}
\end{equation}
To identify the spaces defined by (completion with respect to) the seminorms
$\norm{|\nabla f|}_p$ and $\norm{\Delta^{1/2} f}_p$  on $\ccc_0^\infty(M)$ for some $p\in
(1,\infty)$, it is enough to prove that there exist $0<c_p\le
C_p<\infty$ such that for all $f\in
\ccc_0^\infty(M)$
\begin{equation}\label{comparison}
c_p\norm{\Delta^{1/2} f}_p\le \norm{|\nabla f|}_p\leq C_p\norm{\Delta^{1/2} f}_p. 
\end{equation}
Note that the right-hand inequality may be reformulated by saying that 
 the Riesz transform $\nabla \Delta^{-1/2}$
   is bounded from $L^p(M,\mu)$ to the   space of  $L^p$ vector fields\footnote{In the case where $M$ has finite measure,
 instead of $L^p(M)$, one
has to consider the space
$L^p_0(M)$ of  functions in
$L^p(M)$ with mean zero; this modification will be implicit in what follows.},   in other
words
$$\||\nabla \Delta^{-1/2}f|\|_p \le C_p
 \|f\|_p 
\eqno(R_p)$$
for some constant $C_p$ and all $f\in
\ccc_0^\infty(M)
$.
Note that  $(R_2)$ is trivial from (\ref{triv}).
It is well-known (see \cite{Ba0}, Section 4,
or
\cite{CDfull}, Section 2.1) that the right-hand inequality in \eqref{comparison}  
 implies the reverse inequality
$$
\norm{\Delta^{1/2} f}_q\leq C_p\norm{|\nabla f|}_q,
$$
for all $f\in
\ccc_0^\infty(M)
$, where $q$ is the conjugate exponent of $p$.
Hence, \eqref{comparison} for all $p$ with $1<p<\infty$ follows from $(R_p)$  for all $p$ with $1<p<\infty$.
More generally,    if $(R_p)$ holds  for
$1<p<p_0$ (with $2<p_0\le\infty$), one obtains the equivalence \eqref{comparison} and the 
identification of first order Sobolev spaces for
$q_0<p<p_0$,
$q_0$ being the conjugate exponent to $p_0$.  

Under local assumptions on the manifold, one can hope for the inhomogeneous analog of the
equivalence \eqref{comparison}, namely
\begin{equation}\label{comparisonloc}
c_p\left(\norm{\Delta^{1/2}f}_p+\norm{f}_p\right)\le \norm{|\nabla f|}_p+\norm{f}_p\leq
C_p\left(\norm{\Delta^{1/2}f}_p+\norm{f}_p\right),
\end{equation}
for all $f\in
\ccc_0^\infty(M)
$. It suffices then to study the boundedness on $L^p$ of the local Riesz transform
$\nabla(\Delta+a)^{-1/2}$ for some $a>0$ large enough. If, in addition, the bottom of the
spectrum of the Laplace-Beltrami is positive, that is, if
\begin{equation}\label{sp}
\||\nabla f|\|_2 \ge  \lambda \|f\|_2
\end{equation} for some positive
real number
$\lambda$ and  all $f\in
\ccc_0^\infty(M)
$, one can recover 
 \eqref{comparison}  from \eqref{comparisonloc} (see  \cite{CD}, p.1154).

\subsection{Main results}\label{results}

Let us first recall the result of \cite{CD} which deals with $(R_p)$ for $1<p<2$.
Denote by  $B(x,r)$ the open ball of radius $r>0$ and center $x\in M$, and by $V(x,r)$ its
measure $\mu(B(x,r))$. One says  that $M$ satisfies the doubling
property if for all $x\in M$ and $r>0$
$$
V(x,2r)\le C\,V(x,r).\eqno (D)
$$
Denote by  $p_t(x,y)$, $t>0$, $x,y\in M$, the heat kernel of $M$, that is the kernel of  the heat semigroup $e^{-t\Delta}$. 
One says
 that
$M$  satisfies  the on-diagonal heat kernel upper estimate   if 
$$p_t(x,x)\le \frac{C}{ V(x,\sqrt{t})}, \eqno(DU\!E)$$
for all $x\in M$, $t>0$ and some constant $C>0$.

\begin{theo}\label{CD} Let $M$ be a complete non-compact Riemannian manifold. Assume that $(D)$ and
$(DU\!E)$ hold. Then  the Riesz transform is bounded on $L^p$ for $1<p<2$.
\end{theo}

It is also shown in
\cite{CD} that  the Riesz transform is unbounded on $L^p$  for every $p>2$ on the manifold consisting of two copies of the
Euclidean plane glued smoothly along their unit circles, although it satisfies
$(D)$ and $(DU\!E)$. A stronger
assumption is therefore required in order to extend Theorem \ref{CD} to the range $p>2$.

It is well-known (\cite{G2}, Theorem 1.1) that, under $(D)$,  $(DU\!E)$ self-improves into the off-diagonal upper estimate: 
$$p_t(x,y)\le \frac{C}{ V(
y,\sqrt{t})}\exp\left(-c  \frac{d^2(x,y)}{
t}\right), \eqno (U\!E)
$$
for all $x,y\in M$, $t>0$ and some constants $C,c>0$.
A natural  way to strengthen the assumption is to impose a lower bound of the same size, that is
the full Li-Yau type estimates
$$
 \frac{c}{ V(y,\sqrt{t})}\exp\left(-C  \frac{d^2(x,y)}{
t}\right)\le p_t(x,y)\le \frac{C}{ V(
y,\sqrt{t})}\exp\left(-c  \frac{d^2(x,y)}{
t}\right), \eqno (LY)
$$
for all $x,y\in M$, $t>0$ and some constants $C,c>0$. It is known since  \cite{LY} that such estimates hold on manifolds with
non-negative Ricci curvature. Later, it has been proved in \cite{S2} that
$(LY)$ is equivalent to
the conjunction of $(D)$ and the Poincar\'e inequalities $(P)$ which we introduce next.
 
We say that $M$ satisfies the (scaled) Poincar\'e inequalities if there exists $C>0$ such that, for every ball $B=B(x,r)$ and
every
$f$ with
$f,
\nabla f$ locally square integrable,
 $$
\int_{B}|f-f_{B}|^2\,d\mu\le Cr^2\int_{B}
|\nabla f|^2\,d\mu, \eqno (P)
$$
where $f_E$ denotes the mean of $f$ on $E$. 

However, it follows from the results in  \cite{LHQ} and \cite{CouLi} that  even $(D)$ and $(P)$  do
not suffice for the Riesz transform to be bounded on $L^p$ for any $p>2$. 

In fact, there is also an easy necessary condition for $(R_p)$ to hold. Indeed, $(R_p)$ implies
$$\||\nabla e^{- t\Delta}f|\|_p \le C_p\|\Delta^{1/2}e^{- t \Delta}f\|_p \le
\frac{{C_p}^\prime}{\sqrt t}
\|f\|_p, $$
for all $t>0$, $f\in L^p(M,\mu)$, since, according to \cite{topics},  the heat semigroup is analytic on $L^p(M,\mu)$. And this
estimate may not hold for $p>2$, even in presence of $(D)$ and $(P)$.  

Our main result is that,  under $(D)$ and $(P)$, this condition is sufficient for $(R_q)$, $2<q<p$.
Denote by $\|T\|_{p\to p}$ the norm of a bounded sublinear operator $T$
from  $L^p(M,\mu)$ to itself. 
 
\begin{theo}\label{mainp>2} Let $M$ be a complete non-compact Riemannian manifold satisfying $(D)$
and
$(P)$ $($or, equivalently, $(LY))$. If for some  $p_0 \in (2,\infty]$ there
exists $C>0$ such that, for all $t>0$, 
 $$\||\nabla e^{- t\Delta}|\|_{p_0\to p_0} \le  \frac{C}{\sqrt t} , \eqno(G_{p_0})
$$
then  the Riesz transform  is bounded on $L^p$ for $2<p<p_0$.
\end{theo}

We therefore obtain the announced  necessary and sufficient 
condition as follows. 

\begin{theo}\label{nsc} Let $M$ be a complete non-compact Riemannian manifold satisfying $(D)$
and
$(P)$ $($or, equivalently, $(LY))$. Let $p_0 \in (2,\infty]$. The following assertions are
equivalent:
\begin{enumerate}
\item
  For all  $p \in (2,p_0)$, there exists $C_p$ such that  
 $$\||\nabla e^{- t\Delta}|\|_{p\to p} \le  \frac{C}{\sqrt t},
$$
for  all $t>0$
$($in other words, $(G_p)$ holds for all $p \in (2,p_0))$.
\item The Riesz transform $\nabla \Delta^{-1/2}$ is bounded on $L^p$ for $p\in (2,p_0)$. 
\end{enumerate}
\end{theo}

Notice that we do not draw a conclusion
for $p=p_0$.  It has been shown in \cite{LHQ}  that there exist (singular) manifolds, namely conical  manifolds  with closed
basis, such that  $(R_p)$ hold  if  $1< p< p_0$,  but not for $p\ge p_0$, for some $p_0\in(2,\infty)$, and it has been observed
in \cite{CouLi} that these manifolds do satisfy $(D)$ and $(P)$; in that case, $(G_{p_0})$ does not hold either.
Strictly speaking, these manifolds are not complete, since they have a point singularity. 
But one may observe that our proofs do not use completeness in itself, but rather stochastic completeness,
that is the property
\begin{equation}
\int_Mp_t(x,y)\,d\mu(y)=1,\ \forall\,x\in M,\,t>0,\label{stoco}
\end{equation}
which does hold for complete manifolds satisfying $(D)$, or more generally condition $(E)$ below  (see \cite{Grigostoc}),
but also for  conical  manifolds  with closed basis.

\bigskip

    It follows  from Theorems
\ref{nsc} and 
\ref{CD} that  the assumptions of Theorem
\ref{nsc} are sufficient for $(R_p)$ to hold for $p\in(1,p_0)$, and for  the equivalence \eqref{comparison}  to hold for
$p\in(1,p_0)$ where
$q_0$ is the conjugate exponent to $p_0$. In the case $p_0=\infty$, one can formulate  a sufficient condition for
$(R_p)$ and \eqref{comparison} in the full range
$1<p<\infty$ in terms of pointwise bounds of the heat kernel and its gradient.

\begin{theo}\label{maincor} Let $M$ be a complete non-compact Riemannian manifold satisfying $(D)$
and
$(DU\!E)$. If there exists $C$ such that, for all $x,y\in M$, $t>0$, 
$$|\nabla_x\,p_t(x,y)|\le  \frac{C}{
\sqrt{t}\left[ V(y,\sqrt{t})\right]}, \eqno(G)$$ then the Riesz transform is bounded on $L^p$ 
 and the equivalence \eqref{comparison} holds for
$1<p<\infty$.
\end{theo}

We have seen that under $(D)$,   $(DU\!E)$ implies $(U\!E)$, which, together  with $(G)$, implies 
 the full estimate $(LY)$ (see for instance \cite{LY}).
We shall see in Section \ref{other} that $(G)$ implies 
$(G_p)$ for all $p\in(2,\infty)$. This is why Theorem \ref{maincor} is a corollary of  Theorems
\ref{nsc} and 
\ref{CD}.

\bigskip

 Our results admit local versions.
We  say that $M$ satisfies the exponential  growth  property $(E)$
if for all $r_0>0$, for all $x\in M,\, \theta>1,\, r\le r_0$,
$$
V(x,\theta r)\leq m(\theta) V(x,r), 
\eqno(E)
$$
where  $m(\theta)= C e^{c\theta}$ for some $C\ge 0$ and $c> 0$ depending on
$r_0$.
Note that this implies the local doubling property $(D_{loc})$: for all $r_0>0$ there
exists $C_{r_0}$ such that  for all
$x\in M,\, r\in ]0,r_0[$,
$$
V(x,2r)\leq C_{r_0}\,V(x,r).
\eqno(D_{loc})
$$
We write $(DU\!E_{loc})$ for the property  $(DU\!E)$ restricted to small times (say,
$t\le 1$).

We  say that $M$ satisfies the local Poincar\'e property $(P_{loc})$ if for all $r_0>0$ there
exists $C_{r_0}$ such that for every ball $B$ with radius $r\le r_0$ and every  function 
$f$ with $f,\nabla f$  square integrable on $B$,
$$
 \int_{B}|f-f_{B}|^2\,d\mu\le C_{r_0}r^2\int_{B}
|\nabla f|^2\,d\mu. \eqno (P_{loc})
$$

\begin{theo}\label{mainp>2loc} Let $M$ be a complete non-compact Riemannian manifold satisfying 
$(E)$ and
$(P_{loc})$.  If for some 
$p_0
\in (2,\infty]$ and $\alpha\ge 0$, and for all $t>0$,  
$$\||\nabla e^{- t\Delta}|\|_{p_0\to p_0} \le  \frac{Ce^{\alpha t} }{\sqrt t}, \eqno(G_{p_0}^{loc})
$$ then 
  the local Riesz transform $\nabla (\Delta +a)^{-1/2}$ is bounded on $L^p$ for $2<p<p_0$ and 
$a>\alpha$.
\end{theo}

As a consequence, we can state the

\begin{theo}\label{nscloc} Let $M$ be a complete non-compact Riemannian manifold satisfying 
$(E)$ and
$(P_{loc})$.  Let $p_0
\in (2,\infty]$. Then the following assertions are equivalent:
\begin{enumerate} 
\item For  all $p
\in (2,p_0)$, all $t>0$ and some $\alpha\ge 0$ 
$$\||\nabla e^{- t\Delta}|\|_{p\to p} \le  \frac{C_pe^{\alpha t} }{\sqrt t}.
$$ 
\item
  The local  Riesz transform $\nabla (\Delta +a)^{-1/2}$ is bounded on $L^p$ for $2<p<p_0$
and  some $a>0$.
\end{enumerate}
\end{theo}

Taking into account the local result in \cite{CD}, and denoting by $(G_{loc})$ condition $(G)$ restricted to small times, the
main corollary for the full range $1<p<\infty$  is 

\begin{theo}\label{maincorloc} Let $M$ be a complete non-compact Riemannian manifold satisfying
$(E)$,
$(DU\!E_{loc})$ and  $(G_{loc})$. Then, for $a>0$ large enough, the local  Riesz
transform
$\nabla(\Delta+a)^{-1/2}$ is bounded on
$L^p$   and the equivalence \eqref{comparisonloc} holds  for
$1<p<\infty$.
\end{theo}

Finally, thanks to the argument in \cite{CD}, p.1154, one obtains

\begin{theo}\label{pangap0} Let $M$ be a complete non-compact Riemannian manifold satisfying 
$(E)$,
$(P_{loc})$ 
 and \eqref{sp}. Assume that $(G_{p_0}^{loc})$ holds for some $p_0\in (2,\infty]$.
Then $(R_p)$ holds for all $p\in (1,p_0)$, and \eqref{comparison} holds for all $p\in (q_0,p_0)$, where $q_0$ is the conjugate
exponent to $p_0$.
\end{theo}

\noindent and, in particular,

\begin{theo}\label{pangap} Let $M$ be a complete non-compact Riemannian manifold satisfying 
$(E)$,
$(DU\!E_{loc})$,  $(G_{loc})$ 
 and \eqref{sp}.
Then \eqref{comparison} holds for all  $p$, $1<p<\infty$.
\end{theo}

The core of this paper is concerned with the proof of Theorems \ref{mainp>2} and 
\ref{maincor} and of their local versions
Theorems \ref{mainp>2loc} and \ref{maincorloc}. Before going into details, we comment  on anterior results, on the nature
of our assumptions, and  on our
method.

\subsection{Anterior results} \label{history}

The state of the art  consists so far of a list
of (quite interesting and typical) examples
  with \textit{ad hoc} proofs  rather than a general theory.
These examples essentially  fall into three categories:

\begin{itemize}

\item[I.]  Global statements for manifolds with at most polynomial growth

\subitem{1.} manifolds with non-negative Ricci curvature (\cite{Ba}, 
\cite{Basurv}).

\subitem{2.} Lie groups with
polynomial volume growth endowed with a sublaplacian (\cite{A}).

\subitem{3.} co-compact covering manifolds with
polynomial growth deck transformation group (\cite{ND}).

\subitem{4.} conical manifolds with compact basis without boundary
(\cite{LHQ}).

\item[II.] A local statement

\subitem{5.} manifolds with  Ricci curvature bounded below (\cite{Ba}, 
\cite{Basurv}).
\medskip

\item[III.]  Global statements for manifolds where the bottom of the 
spectrum is positive

\subitem{6.}  Cartan-Hadamard manifolds where the Laplace operator is 
strictly positive,  plus   bounds on the curvature  tensor and
its two first derivatives (\cite{Lo}).

\subitem{7.} unimodular, non-amenable Lie
groups (\cite{Lo1}). 

\end{itemize}

Note that, for results concerning Lie groups in the above list,  one can  consider not only the case where they are
endowed with a translation-invariant Riemannian metric, but also the case where  they are endowed with a  sublaplacian,
that is a sum of squares of invariant vector fields satisfying the H\"ormander condition. For more on this,
see for instance \cite{A}. Although we did not introduce this framework, for the sake of brevity, our proofs do work without
modification in this setting also, as well as, more generally, on a manifold endowed with a subelliptic sum of squares of
vector fields. 

In cases I and III, the conclusion is the boundedness of the Riesz transform, hence the  seminorms equivalence 
(\ref{comparison}), for all
$p\in(1,\infty)$.

In case II,
 the conclusion involves local Riesz transforms, or the equivalence (\ref{comparisonloc}) of inhomogeneous Sobolev norms.
Note that an important feature of Bakry's result in this case (say, \cite{Ba}, Theorem 4.1) is the weakness of the assumption:
neither  positivity of the injectivity radius nor bounds on the derivatives of the curvature tensor are assumed. By contrast,
for bounded geometry manifolds, (\ref{comparisonloc}) follows  easily from the Euclidean result by patching. 

We feel that there is a logical order between I, II, III:
the results in II are nothing but  local versions of I, and III 
follows easily from II if one uses the additional assumption
on the spectrum of the Laplace operator. One may observe that the 
above results  were in fact obtained in a quite different
chronological order.

The results in I are covered by Theorem \ref{maincor}, the one in II by Theorem \ref{maincorloc},  and the ones in III by 
Theorem \ref{pangap}.
Let us  explain now in each of the above
situations where the required assumptions come from.

In case 1, the doubling property follows from Bishop-Gromov volume comparison (\cite{Chavel}, Theorem 3.10),
and the heat kernel bounds including $(G)$ from \cite{LY}.
Note  however that an important additional outcome of Bakry's method in \cite{Ba}, 
\cite{Basurv} (see also \cite{Ba0} for a more abstract setting) is the independence  of constants with respect to the
dimension. See the comments in Section \ref{rela}.

In case 2, the doubling property is obvious, the heat kernel  upper 
bound follows from
\cite{V} (but one has nowadays much simpler proofs, see for instance 
\cite{CGP} for an
exposition), and the gradient bound from
\cite{saloffpoly}.  The proof
in \cite{A} is much more complicated than ours; it requires some 
structure   theory of Lie groups
as well as  a  substantial amount of homogeneization theory.

In case 3,   $(D)$ is   again obvious  since such a manifold has polynomial 
volume growth, $(DU\!E)$ is well-known (it
can be extracted from the work of Varopoulos, see for instance
\cite{Varcov}, but nowadays one can write down a simpler proof by 
using \cite{CDisc} or \cite{CS}; see for instance
\cite{G3}, Theorem 7.12)
and
$(G)$ is proved in
\cite{ND} ({Added after acceptation:} another simpler proof of $(G)$ has been proposed recently in \cite{ND3}.)  The boundedness of the Riesz transform is directly 
deduced in \cite{ND} by using further specific properties  of this
situation.

In case 4, the boundedness of the Riesz transform is obtained for a range $(1,p_0)$ of values 
of $p$, and is shown to be false outside this range. It follows from
\cite{LHQ2} that $(G_p)$ holds for $1<p<p_0$, hence yielding
with our result a simple proof of the main results in \cite{LHQ}. 
By direct estimates from below on $\nabla 
e^{-t\Delta}$ as in \cite{LHQ2}, one can also recover the negative
results   for $p\ge p_0$ (see \cite{CouLi}).

Cases 6, 7 are covered by Theorem \ref{pangap}. In case 6, we get 
rid of specific regularity  assumptions on the
curvature tensor.  As far as case 7 is concerned, for more recent  results related to Lie groups with exponential growth,  see
\cite{LM1},
\cite{LM2}, \cite{GS}, \cite{HS}. Note that the groups considered there are either non-unimodular or non-amenable, which 
allows  
reduction to a local problem by
use of the positivity of the bottom of the spectrum.  

Let us finally mention a few results  which are {\em not} covered by our methods: in \cite{LHQ}, conical  manifolds with compact basis with boundary is considered; in that case, the conical manifold is not  complete.
The case where the basis is non-compact has  been considered in \cite{LHQ}, and studied further in \cite{LL}; here the
volume of balls with finite radius may even be infinite.   In \cite{LHQ3}, the $L^p$ boundedness of the Riesz transform for all
$p\in (1,\infty)$ is obtained for a specific class of
manifolds with exponential volume growth, namely cuspidal  manifolds
with compact basis without boundary. In \cite{LXD}, Theorem 2.4, the boundedness of Riesz transform for
$p>2$ is  proved for a class of Riemannian manifolds with a certain amount of  negative curvature; here doubling is not
assumed, and the main tool is Littlewood-Paley theory, as in \cite{CDfull}.

\subsection{About our assumptions}\label{assumption}

We discuss here the meaning and relevance of our assumptions.

Let us begin with the basic assumptions on the heat kernel.
The  two  assumptions  $(D)$ and $(DU\!E)$ in Theorem \ref{CD}
  are   equivalent, according to \cite{G}, to the
so-called relative Faber-Krahn inequality
$$\lambda_1(\Omega)\ge
\frac{c}{r^2}\left(\frac{V(x,r)}{\mu(\Omega)}\right)^{2/\nu}, \eqno(FK)$$
for  some $c,\nu>0$, all $x\in
M$, $r>0$, $\Omega$ smooth  subset of $B(x,r)$. Here 
$\lambda_1(\Omega)$ is the first eigenvalue of the Laplace operator  on
$\Omega$ with Dirichlet boundary conditions :
$$\lambda_1(\Omega)=\inf\left\{\frac{\int_\Omega|\nabla u|^2}{\int_\Omega
u^2},\ u\in
C^\infty_0(\Omega)\right\}.
$$
In the sequel, we shall sometimes denote by $(FK)$ the conjunction of $(D)$ and  $(DU\!E)$.
Also, we have recalled that the conjunction of $(D)$ and $(P)$ is equivalent to $(LY)$. 
It is worthwhile to note that, contrary to the non-negativity of the Ricci curvature, $(D)$ and $(P)$ are  invariant under
quasi-isometry, which is not obvious  to check directly on $(LY)$.

We may question the 
relevance of this group of assumptions to Riesz transform bounds; as a matter of fact, $(G_p)$ is needed but neither
$(FK)$ for $p<2$ in
\cite{CD} nor $(D)$ and $(P)$ for $p>2$ are known to be necessary. However, it seems out of reach as of today to prove such
bounds without some minimal information on the heat kernel. These assumptions are reasonable
 for the moment but we think they can be weakened.  One direction 
  is to replace the doubling condition by  exponential volume growth (without positivity of the bottom of
the spectrum). This would mean extending the  Calder\'on-Zygmund theory to the exponential growth
realm. A very promising tentative in this direction is in \cite{HS}, 
although the method has been so far only applied  to a (typical) class of Lie groups having a positive bottom of the
spectrum. Another direction is to pursue
\cite{CDfull} by using Littlewood-Paley-Stein functionals and prove the conjecture stated there.
Again, see the comments in Section \ref{rela}. 
 
Concerning $(P)$, a  minor improvement of our assumptions is that, under $(FK)$ and $(G_{p_0})$,
it may certainly be relaxed to
$L^r$ Poincar\'e inequalities for $r$ large enough so as to guarantee some control on the oscillation of the heat
kernel. We have not tried to go into this direction here. On the other hand, if the manifold has polynomial volume growth, then
as soon as  such weak Poincar\'e  inequalities hold, $(P)$ is necessary for the Riesz transform to be bounded on $L^p$ for $p$
larger than the volume growth exponent
(see \cite{CMorrey}, Section 5).  
\bigskip

We continue with   estimates on the gradient of the heat kernel.

First, we can reformulate the necessary 	and sufficient condition in Theorem \ref{nsc} in terms of integral bounds on the
gradient of the heat kernel, thanks to the following proposition  proved  in Section \ref{other}. 

\begin{pro} \label{proother} Assume that $M$ is a complete non-compact Riemannian manifold satisfying
$(FK)$. Let
$2<p_0\le
\infty$. The following assertions are equivalent
\begin{enumerate}
\item $(G_p)$ holds for all $2<p<p_0$.

\item For all $2<p<p_0$, for all $y\in M$ and $t>0$
\begin{equation}\label{gradlp}
\norm{|\nabla_x\,p_t(.,y)|}_p\le  \frac{C_p}{ \sqrt{t}\left[ V(y,\sqrt{t})\right]^{1-\frac{1}{
p}}}.
\end{equation}
\end{enumerate}
 \end{pro}

As one can see,  the case $p=\infty$ is excluded from this statement. When $p=\infty$,
\eqref{gradlp} is precisely $(G)$ in Theorem \ref{maincor} on which we concentrate now.   
Consider two other  conditions on $\nabla_x\, p_t(x,y)$:
\begin{equation}\label{intgradpt}
\sup_{t>0,\ x\in M\ }\sqrt t \int_M \, |\nabla_x\, p_t(x,y)|\, d\mu(y) <\infty. 
\end{equation}
\begin{equation} \label{gradgub}
|\nabla_x\, p_t(x,y)|\le \frac{C}{\sqrt{t}\,V(y,\sqrt{t})}\exp{\left(-c\frac{d^2(x,y)}{t}\right)},
\end{equation}
 for all $t>0$, 
$x,y\in M$.

It is easy to show that, under $(FK)$,  \eqref{gradgub} $\Longrightarrow$ 
\eqref{intgradpt}  $\Longrightarrow  (G)$.
Indeed \eqref{gradgub} $\Longrightarrow$ 
\eqref{intgradpt} is immediate by integration using $(D)$. Let us note that \eqref{intgradpt} 
is equivalent to 
$$\||\nabla e^{- t\Delta}|\|_{\infty\to\infty} \le  \frac{C }{\sqrt t}.
\eqno(G_\infty)
$$
Note in passing that, by interpolation with $(G_2)$, $(G_\infty)$ implies $(G_p)$ for all $p\in (2,\infty)$.  Then,
\eqref{intgradpt} 
$\Longrightarrow  (G)$ follows by writing
$$
\nabla_x\,p_t(x,y) = \int_M \nabla_x\, p_{t/2}(x,z) p_{t/2}(z,y)\, d\mu(z)
$$
and by direct estimates using $(U\!E)$. One can see that in fact, under $(D)$, \eqref{intgradpt} is equivalent to  
\eqref{gradgub}, but this is another story (see \cite{CSiko}).  

\bigskip

 Next, it is  interesting to observe that the size estimates
$(LY)$ do include already some regularity estimates for the heat kernel, and that $(G)$ is nothing but a slightly stronger form
of this regularity.  More precisely, the estimates
$(LY)$ are  equivalent to a so-called uniform parabolic Harnack principle
 (see \cite{parma}) and, by the same token, they imply 
\begin{equation}
|p_t(x,y)-p_t(z,y)|\le \left(\frac{d(x,z)}{ \sqrt{t}}\right)^\alpha 
\frac{C}{V(y,\sqrt{t})},\label{po}
\end{equation}
 for some $C,c>0$, $\alpha\in (0,1)$, and  all $x, y,z \in M, \,t>0$. 
The additional assumption $(G)$ is nothing but the limit case $\alpha=1$ of 
(\ref{po}). One can therefore summarize the situation by saying that  the H\"older regularity of the heat kernel,
yielded by the uniform parabolic Harnack principle, is not enough in general for the Riesz transform to be bounded on all $L^p$
spaces, whereas Lipschitz regularity does suffice.

Unfortunately, $(G)$  does not have  such a nice geometric characterization as
$(D)$ and
$(P)$. In fact, it is unlikely one can find a geometric description of $(G)$ that is invariant under quasi-isometry. It may 
however be the case that
$(G)$, and even more $(G_p)$, are stable under some kind of perturbation of the manifold, and this certainly deserves
investigation.

Let us make a digression.   We already observed that $(G)$ and $(FK)$ imply
 the full $(LY)$ estimates, hence, under $(FK)$,  \eqref{gradgub}  is equivalent to
\begin{equation} \label{gradptpt}
|\nabla_x\,p_t(x,y)|\le \frac{C}{ \sqrt{t}}\,p_{C't}(x,y)
\end{equation} 
for some constants $C,C'>0$. 
Note that, in the case where $C'=1$, this can be reformulated as 
\begin{equation}\label{logpt}
|\nabla_x\, \log p_t(x,y)|\le \frac{C}{ \sqrt{t}},
\end{equation}
which is one of the fundamental bounds for manifolds with non-negative Ricci curvature (see \cite{LY}).

\bigskip

 Known methods to
prove pointwise gradient estimates   include the Li-Yau method (\cite{LY}, see also \cite{Q} for
generalizations), as well as coupling (\cite{C}), and other probabilistic methods (see for
instance \cite{Pj}) including  the derivation of  Bismut type formulae  which enable one to estimate the logarithmic derivative
 of the heat kernel as in (\ref{logpt}) (see for instance  \cite{EL}, \cite{ST}, \cite{TW1} 
and references therein). Unless
one assumes non-negativity of the curvature, all these methods are limited  so far to small time, more precisely they yield
the crucial factor $\frac{1}{\sqrt{t}}$ only for small time. One may
wonder which large scale geometric features, more stable and less specific than non-negativity of the Ricci curvature, 
would be sufficient to ensure a large time  version of such estimates. 
 A nice statement is that if $M$ satisfies $(FK)$ and if for all $t>0$, $x,y \in M$,
 $$|\nabla_x\,p_t(x,y)|\le C |\nabla_y\,
p_t(x,y)|$$ 
then \eqref{gradgub} (therefore $(G)$)
holds (\cite{G1}, Theorem 1.3, see also the first remark after Lemma \ref{Steve} in Section
\ref{sectionsteve}). Another interesting approach is in \cite{ND}, where $(G)$ is deduced from a discrete regularity estimate,
but here a group invariance is used in a crucial way; see also related results in \cite{I}. The question of finding weaker 
sufficient conditions for   the integrated estimates
$(G_p)$ is so far completely open.  
\bigskip

We note that our work is not the first example of 
the phenomenon that higher integrability of gradients of solutions 
is related to the $L^p$ boundedness of singular integrals. ÊIndeed, 
the property that the gradient of the heat kernel 
satisfies an $L^p$ bound can be thought of 
as analogous to the estimate of Norman Meyers \cite{Ms} concerning the 
higher integrability of gradients of solutions, which in turn is 
connected to Caccioppoli inequalities and reverse H\"older inequalities. 
It has been pointed out by T. Iwaniec in \cite{Iw}  that $L^p$ reverse H\"older/Caccioppoli inequalities for solutions to a
divergence form  elliptic equation $Lu=f$ are equivalent (at least up to endpoints) 
to the $L^p$ boundedness of the 
Hodge projector $\nabla L^{-1} {\rm div}$.\footnote{\,{\bf Added after acceptation:} in fact,   a recent work by Z. Shen \cite{Shen} shows that reverse H\"older inequalities are equivalent to $L^p$ boundedness of the Hodge projector and also to the $L^p$ boundedness of the Riesz transform $\nabla L^{-1/2}$ when $p>2$ and $L$ is a real symmetric uniformly elliptic operator $-div(A \nabla)$ on Lipschitz domains of $\mathbb{R}^n$.    For this, he states a general theorem akin to  our Theorem 2.1 and attributes the method of proof to ideas of Caffarelli-Peral \cite{CP}. In a subsequent paper \cite{AC}, the  two first-named authors of the present article will extend these ideas to the manifold setting and prove actually that  there always is some $p_{0}>2$ for which our  Theorem  1.3 applies.}
ÊOur result, which says that 
the $L^p$ boundedness of the gradient of the heat semigroup 
is equivalent (again up to endpoints), to the $L^p$ boundedness of the 
Riesz transform, is thus in the same spirit. 
\bigskip

Let us finally connect \eqref{gradgub} with properties of the heat kernel on 1-forms. In \cite{CD2},
\cite{CDfull},   the boundedness of the Riesz transform on $L^p$ is proved for $2<p<\infty$ under
$(FK)$ and the assumption that the heat kernel on $1$-forms is dominated by the heat kernel on
functions: for all $t>0,\,\omega\in
\ccc^\infty T^*M$,
$$|e^{-t\overrightarrow{\Delta}}\omega|\le C e^{-ct\Delta}|\omega|,
$$
(the case $C=c=1$ of this estimate corresponds to non-negative Ricci curvature). 
It would certainly be  interesting  to investigate the class of manifolds where this domination
condition holds;    unfortunately, this  is in general  too
strong a requirement, since for instance it does not hold  for nilpotent Lie groups, as is shown in \cite{Rm},
\cite{Rm1}, whereas  on such groups the Riesz transforms are known to be bounded on $L^p$ for all $p\in (1,\infty)$ (\cite{A}).

It is conjectured in \cite{CDfull} that the same result is true under a weak 
commutation between the gradient and the semigroup (since commutation is too much to ask, of course)
that is the restriction of the above domination condition to {\em exact} forms: for all $t>0$,  $f\in
\ccc_0^\infty(M)$,
$$|\nabla e^{-t{\Delta}}f|\le C e^{-ct\Delta}|\nabla f|,
$$
and even its weaker, but more natural, $L^2$ version:
\begin{equation}
|\nabla e^{-t{\Delta}}f|^2\le C e^{-ct\Delta}(|\nabla f|^2).\label{conj}\end{equation}
({Added after acceptation:} we mention a paper by 
Driver and Melcher \cite{DM} where the $L^p$ versions of such inequalities (with $c=1$)
are proved on the Heisenberg group ${\mathbb H}^1$ for all $p>1$ by probabilistic methods.) This is what we prove here, in the class of
manifolds with
$(FK)$,   as a consequence of Theorem \ref{maincor}, since then \eqref{conj} is equivalent to 
\eqref{gradgub} as we show in Lemma \ref{Steve}, Section \ref{sectionsteve}, when we give a 
simpler argument for proving Theorem \ref{maincor} with \eqref{gradgub} instead of $(G)$.

\subsection{About our method}\label{method}

Let us emphasize several features of our method.

First, we develop an appropriate machinery to treat
operators beyond the classical Calder\'on-Zygmund operators. Indeed, our operators no longer have
H\"older continuous kernels, as this is often too demanding in applications:
 the kernel of the Riesz transform is formally given by $$\int_0^\infty \nabla_x\, p_t(x,y)\frac{dt}{\sqrt{t}},$$ and 
condition $(G)$
 is just an upper bound which does not require  H\"older regularity on $\nabla_x\, p_t(x,y)$ in a
spatial variable. 
Geometrically, this makes a big difference since pointwise upper bounds on the oscillation in $x$ of $\nabla_x\,p_t(x,y)$ seem
fairly unrealistic for large time (see
\cite{LJ}). The loss of H\"older continuity is compensated by  a built-in regularity
property from the semigroup
$e^{-t\Delta}$.  Such an idea, which originates 
from 
\cite{He},  has been formalised in \cite{DuMcIn} for boundedness results in the range $1<p<2$ and
is actually used in
\cite{CD} to derive Theorem \ref{CD}.  However, this method does not apply to our situation as
$p>2$; a duality  argument would not help us either as we would have to make assumptions on the
semigroup acting on 1-forms  as explained in Section \ref{assumption}; this would bring us back to the state of the art in
\cite{CDfull} (see
\cite{Siko} for this approach to the results in \cite{CDfull}). But recently, it was shown in   
\cite{M} that this regularity property can be used for $L^p$
results in the range
$p>2$ by employing  good-$\lambda$ inequalities
as in 
\cite{FS} for an \textit{ad hoc} sharp maximal function; this may be seen as the basis to an $L^\infty$ to $BMO$ version of the
$L^1$ to weak $L^1$ theory in \cite{DuMcIn}. Note that, in this connection, the usual BMO theory requires too strong
assumptions and can only work in very special situations (see
\cite{Ch}).

Second, our method works for the  usual full range
$2<p<\infty$ of values of $p$ and also for  a limited range  $2<p<p_0$.  This is important as in
applications (to Riesz transforms on manifolds or to  other situations, see \cite{Au})  the
operators may no longer have  kernels with pointwise bounds! This is akin to results  recently
obtained in
\cite{BK1} for
$p<2$ and non-integral operators,  which generalize
\cite{DuMcIn}; in this circle of ideas, see also \cite{HM}. 
Here, we  state a general theorem valid in arbitrary range of $p$'s above 2, and its local version
(Theorems \ref{lp} and \ref{lploc}, Section \ref{basic}).

Third, we use very little of the
differential structure on manifolds, and in particular we do not use the heat kernel on 1-forms as
 in \cite{Ba}, \cite{Basurv},  or \cite{CDfull}. As a matter of fact, our method is quite general, and enables one to prove
the $L^p$ boundedness of a Riesz transform of the form $\nabla L^{-1/2}$ as soon as the following ingredients are available: 

1. Doubling measure
 
2. Scaled Poincar\'e inequalities
 
3. $e^{-tL} 1 = 1$ 

4. ellipticity, a divergence form structure, and "integration by parts" 
(in other words, the ingredients necessary to prove 
Caccioppoli type inequalities)
 
5. $L^p$ boundedness of $ \sqrt{t} \nabla e^{-tL}$

6. $L^2$ bound for the Riesz transform. 

In particular, the method applies equally well to accretive (i.e. elliptic) divergence form operators on $\mathbb{R}^n,$ in which
case the 
$L^2$ bound is equivalent to the solution of the square root problem 
of Kato \cite{AHLMT}. ÊOne of the present authors (Auscher) will present the details 
of this case (as well as related results) 
in a forthcoming article \cite{Au}. 

The method is also subject to further
extensions to other settings such as general Markov diffusion semigroups on metric measure spaces, or discrete Laplacians on
graphs.   See Section
\ref{rela}.

\section{Singular integrals and a variant of the
sharp maximal
function}\setcounter{equation}{0}\label{basic}

In this section, $(M,d,\mu)$ is a  measured metric space.
We  denote as above by  $B(x,r)$ the open ball of radius $r>0$ and center
$x\in M$,  which we assume to be always  of finite $\mu$-measure. We state and prove 
a criterion for $L^p$ boundedness,  with $p>2$,  for operators such as   singular integrals or 
quadratic expressions.  We also give a local analog of this criterion. 

\subsection{The global criterion}\label{global}

We   say that $M$ satisfies the doubling
property (that is $(M,d,\mu)$ is of homogeneous type in the terminology of 
\cite{CW}) if there exists a constant $C$ such that, for all $x\in M,\,r>0$, 
$$
\mu(B(x,2r))\le C\,\mu(B(x,r)).\eqno (D)
$$

Consider a sublinear operator acting on $L^2(M,\mu)$. We are  going to prove a general statement
 that allows one to obtain a bound for
its operator norm on $L^p(M,\mu)$   for a fixed
$p>2$.  Such techniques originate,  in a Euclidean setting, in \cite{FS}
(see also
\cite{SFS}) by use of the sharp maximal function and good-$\lambda$ inequalities. It is proved in
\cite{M} that,  in the definition of the sharp function,  the average
over balls can  replaced by  more general  averaging operators depending on the context, and that 
the ideas of \cite{FS} can be adapted. Our method  is based on that of \cite{M}.

Denote by  ${\cal M}$  the 
Hardy-Littlewood maximal operator
$${\cal M}f(x)=\sup_{B\ni x}\frac{1}{\mu(B)}\int_B|f|\,d\mu,$$
where $B$ ranges over all open balls containing $x$.

\begin{theo}\label{lp} Let $(M, d,\mu)$ satisfy $(D)$ and let $T$ be a sublinear operator which is bounded on $L^2(M,\mu)$. Let
$p_0\in (2,\infty]$. Let  $A_r$, $r>0$,
be  a family of linear operators acting on $L^2(M,\mu)$.
     Assume   
\begin{equation}
\left(\frac{1}{\mu(B)}\int_B|T(I-A_{r(B)})f|^2\,d\mu\right)^{1/2} \le C\big( {\cal M}
(|f|^2)\big)^{1/2}(x) ,\label{pointwise}
\end{equation}
and 
\begin{equation}
\left(\frac{1}{\mu(B)}\int_B |TA_{r(B)}f|^{p_0} \,d\mu\right)^{1/p_0}\le  C\big({\cal M}
(|Tf|^2)\big)^{1/2} (x),\label{doma}
\end{equation}
for all $f\in L^2(M,\mu)$, all $x\in M$ and all balls $B\ni x$,  $r(B)$ being the radius of $B$.
If  $2<p<p_0$ and  $Tf \in L^p(M,\mu)$ when $f\in L^p(M,\mu)$,  then $T$ is of strong type $(p,p)$ and its
operator norm is bounded by a constant depending only on  its $(2,2)$ norm, on  the constant in $(D)$, on $p$ and $p_0$, and on
the constants in
\eqref{pointwise} and
\eqref{doma}.
\end{theo}

\noindent{\bf Remarks:}

-If $p_0=\infty$, the left-hand side of \eqref{doma} should be understood as the essential supremum
 $\sup_{y\in B}|TA_{r(B)}f(y)|$.

- The operators $A_r$ play the role of approximate identities (as $r\to 0$). Notice that  the
regularized version $TA_r$ of
$T$ is controlled by the maximal function of $|Tf|^2$ which may be surprising at first
sight since
$T$ is the object under study. The improvement from 2 to $p_0$ in the exponents expresses
a regularizing effect of $A_r$.

- Define, for $f\in
L^2(M,\mu)$, 
$${\cal M}^{\#}_{T,A}f(x)=\sup_{B\ni
x}\left(\frac{1}{\mu(B)}\int_B|T(I-A_{r(B)})f|^2\,d\mu\right)^{1/2},$$ where the supremum is taken
over all balls $B$ in $M$ containing $x$,  and $r(B)$ is  the radius of $B$. This is (a variant of) the
substitute to the sharp function alluded to above. Assumption \eqref{pointwise} means that ${\cal M}^{\#}_{T,A}f$ is
controlled pointwise by $({\cal M}
(|f|^2))^{1/2}$. In fact, rather than the exact form of the control, what matters is that ${\cal
M}^{\#}_{T,A}$ is of strong type
$(p,p)$ for the desired values of $p$.

- Note that we assumed that $T$ was already acting on $L^p(M,\mu)$ and then we obtained boundedness and
a bound of its norm. In practice, this theorem is applied to suitable approximations of $T$, the
uniformity of the bound allowing a limiting argument to deduce $L^p$ boundedness of $T$ itself. 

-A careful reader will notice that in the proof below, the $L^2$ bound for $T$ is explicitely used only if 
$M$ has finite volume; but in practice, the verification of the assumptions (\ref{pointwise}) and (\ref{doma})
requires the $L^2$ boundedness of $T$ (and $A_r$) anyway.

\bigskip

Let us now prove two lemmas  inspired from \cite{M} but with modifications to allow a treatment
at a given exponent $p$ and for a right regularization (see the remark after the proof).  The first one
is a so-called good-$\lambda$ inequality. 
For simplicity, we normalize  in the following proof the constants 
in
assumptions \eqref{pointwise} and
\eqref{doma} to one.

\begin{lem}\label{gooda} Let $(M, d,\mu, A_r, p_0)$ and $T$ be as above. 
Assume  that 
${\rm (\ref{doma})}$ holds. There exist $K_0>1$ and $C>0$ only depending on   $p_0$ and the constant in $(D)$,
such that, for every $\lambda>0$, every
$K>K_0$ and $\gamma>0$, for every ball $B_0$ in $M$ and every function $f\in L^2(M,\mu)$ such
that there exists $x_0\in B_0$ with ${\cal M}(|Tf|^2)(x_0)\le \lambda^2$, then 
\begin{equation}
\mu\left(\{x\in B_0; {\cal M}(|Tf|^2)(x)>K^2\lambda^2,\ {\cal M}_{T,A}^\#f(x)\le
\gamma\lambda\}\right)\le
C(\gamma^2 + K^{-p_0})
\,\mu(B_0).
\label{good} 
\end{equation}
\end{lem}

\noindent{\bf Proof:} Let us assume first that $p_0<\infty$. Set $$E=\{x\in B_0; {\cal
M}(|Tf|^2)(x)>K^2\lambda^2,\ {\cal M}_{T,A}^\#f(x)\le
\gamma\lambda\}.$$ 
From (\ref{doma}) and the hypothesis that ${\cal M}(|Tf|^2)(x_0)\le \lambda^2$, one has
$$
\int_{3B_0} |TA_{r_0}f|^{p_0} \,d\mu \le \lambda^{p_0}  \,\mu(3B_0)
$$ where $r_0=r(3B_0)$. 
Denote  
$$\Omega= \{ x \in M; {\cal M}(|TA_{r_0}f|^2\chi_{3B_0})(x)> J^2\lambda^2\}
$$ where $J$ is a positive constant to be chosen. By the weak type $(p_0/2, p_0/2)$ of the maximal
operator, we have
$$
\mu(\Omega) \le  \frac{C}{J^{p_0}\lambda^{p_0}}  \int_{3B_0} |TA_{r_0}f|^{p_0} \,d\mu \le C
J^{-p_0} \mu(3B_0).
$$
Now, we want to estimate $\mu(E\setminus \Omega)$. We remark that by definition, if $x\in E\setminus
\Omega$, then
\begin{equation}
{\cal M}(|TA_{r_0}f|^2\chi_{3B_0})(x) \le J^2\lambda^2.\label{edd}
\end{equation}
We  first  prove that  there exists $c_0$ only depending on $(D)$ such that, if $c_0K^2>1$, then 
for
every $x\in E$,
\begin{equation}
{\cal M}(|Tf|^2\chi_{3B_0})(x)  > K^2\lambda^2.\label{edun}
\end{equation}
Indeed, let $x\in E$. Since ${\cal M}(|Tf|^2)(x)>K^2\lambda^2$, there is a ball $B$ containing $x$ such
that 
\begin{equation}
\int_{B}|Tf|^2\,d\mu >K^2\lambda^2\, \mu(B).\label{pif}
\end{equation}
If $r=2r(B)$, one has $B \subset B(x,r) \subset 3B$, hence for  $c_0$ only depending on the doubling
condition $(D)$ one has $\mu(B) \ge c_0 \mu(B(x,r))$. Therefore,
\begin{equation}
\int_{B(x,r)}|Tf|^2\,d\mu >c_0K^2\lambda^2\, \mu(B(x,r)).\label{pifpaf}
\end{equation}
Since 
 ${\cal M}(|Tf|^2)(x_0)\le \lambda^2$ and $c_0K^2>1$, one can infer that $x_0$ does not
belong to $B(x,r)$.  Therefore $r<2r(B_0)$ and one concludes that $B\subset 3B_0$. 
Together with (\ref{pif}), this yields (\ref{edun}).

Next,  choose $J$ such that $2(J^2+1)=K^2$. 
Then we have, for $x\in E\setminus \Omega$,
\begin{align*}
2(J^2+1)\lambda^2&<{\cal M}(|Tf|^2\chi_{3B_0})(x)\\
&\le 2{\cal
M}(|T(f-A_{r_0}f)|^2\chi_{3B_0})(x)+2{\cal M}(|TA_{r_0}f|^2\chi_{3B_0})(x)\\
&\le 2 {\cal
M}(|T(I-A_{r_0})f|^2\chi_{3B_0})(x)+2J^2\lambda^2,
\end{align*}
and so
$${\cal
M}(|T(I-A_{r_0})f|^2\chi_{3B_0})(x)>\lambda^2.$$ 
Therefore
$$E\setminus \Omega \subset \{x\in M;{\cal
M}(|T(I-A_{r_0})f|^2\chi_{3B_0})(x)>\lambda^2\}.$$
The weak type $(1,1)$ inequality for the Hardy-Littlewood maximal function
 yields
\begin{eqnarray*}
\mu(E\setminus \Omega)&\le& \mu\left(\{x\in M;{\cal
M}(|T(I-A_{r_0})f|^2\chi_{3B_0})(x)>\lambda^2\}\right)\\
&\le&\frac{C}{\lambda^2}\int_M |T(I-A_{r_0})f|^2\chi_{3B_0}\,d\mu\\
&=&\frac{C}{\lambda^2}\int_{3B_0}|T(I-A_{r_0})f|^2\,d\mu\\
&\le&\frac{C}{\lambda^2}\mu(3B_0)\left({\cal M}^{\#}_{T,A}f(x)\right)^2\\
&\le& C\gamma^2 \mu(3B_0).
\end{eqnarray*}
In the last two inequalities, we have used that
$3B_0$ contains $x$, and that
$x
\in E$. Note that $C$ is the weak type $(1,1)$ bound of the
maximal operator and, therefore,  only depends on
$(D)$. 

Altogether, we have obtained that 
$$
\mu(E) \le C(J^{-p_0}+\gamma^2)\mu(3B_0)
$$
provided $K>1$,  $c_0K^2>1$ and $K^2=2(J^2+1)$.
This proves the lemma when $p_0<\infty$.

 If $p_0=\infty$, one deduces from (\ref{doma}) that
$$
 |TA_{r_0}f(x)|  \le {\cal M}(|Tf|^2)^{1/2}(x_0)\le  \lambda  
$$
for  $\mu$-a.e. $x\in 3B_0$. 
Hence
$$
{\cal M}(|TA_{r_0}f|^2\chi_{3B_0})(x) \le \lambda^2
$$
for all $x\in M$, and the set $\Omega$ is empty if $J\ge 1$. The rest of the  proof  proceeds as
before.

\bigskip 

As in \cite{newStein}, Lemma 2, p.152,  the good-$\lambda$ inequality yields comparisons of $L^p$
norms as used in
\cite{M}, Theorem 4.2.

\begin{lem}\label{sharpi} Let $(M, d,\mu, A_r)$ and $T$ be as above. 
Assume  that 
${\rm (\ref{doma})}$ holds.
Then, for  $0< p< p_0$, there exists $C_p$ such that 
\begin{equation}
\|\left({\cal M}(|Tf|^2)\right)^{1/2}\|_p\le C_p\left(\|{\cal M}_{T,A}^\#f\|_p+\|f\|_p\right), \label{sharpa}
\end{equation}
for every $f\in
L^2(M,\mu)$ for which the left-hand side is finite
(if $\mu(M)=\infty$, the term $C_p\|f\|_p$   can be dispensed with in the right-hand side of $(\ref{sharpa})$). 
\end{lem}

\noindent{\bf Proof:} 
 Let $f\in L^2(M,\mu)$ be such that  $\|\left({\cal M}(|Tf|^2)\right)^{1/2}\|_p<\infty$.
For $\lambda>0$, set $$E_\lambda=\{x\in M; {\cal
M}(|Tf|^2)(x)>\lambda^2\}.$$
Set $\lambda_0=0$ if $\mu(M)=\infty$, $\lambda_0=\frac{1}{\mu(M)}\int_M{\cal
M}(|Tf|^2)\,d\mu$ if $\mu(M)<\infty$.
In the latter case, by Kolmogorov inequality (see \cite{Mey}, p.250), $(D)$,  and
the weak type
$(1,1)$ of
${\cal M}$, 
$$
\int_{M}  {\cal M}(|Tf|^2)^{1/2}\, d\mu \le C\mu(M)^{1/2} \||Tf|^2 \|_{1}^{1/2} =
C\mu(M)^{1/2} \|Tf \|_{2}.
$$
Using then the $L^2$ boundedness of $T$, we obtain
$$
\lambda_0 \le \frac{C}{\mu(M)^{1/2}} \|f\|_2 \le  \frac{C'}{\mu(M)^{1/p}}
\|f\|_{p}.
$$
Now fix $K>0$ to be chosen later, and write
$$\|\left({\cal M}(|Tf|^2)\right)^{1/2}\|_{p}^p=I_1+I_2,$$
with
$$I_1=\int_{{\cal
M}(|Tf|^2)
\le K^2\lambda_0^2} \left({\cal M}(|Tf|^2)\right)^{p/2}\, d\mu,$$
$$I_2=\int_{{\cal M}(|Tf|^2)
>K^2\lambda_0^2} \left({\cal M}(|Tf|^2)\right)^{p/2} \, d\mu.$$
Clearly, $I_1$ 
is bounded above by 
$$
K^p\lambda_0^p \mu(M) \le K^p C^p \|f\|_{p}^p
$$
with $C$ depending only on the  constant in $(D)$ and the $L^2$ norm of $T$.
One can treat $I_2$ as follows. 
The Whitney  decomposition (\cite{CW}, Chapter III, Theorem 1.3) for $E_\lambda$
yields, for $\lambda>\lambda_0$, a family of boundedly overlapping balls $B_i$ such that $E_\lambda=\cup_iB_i$. There exists
$c>1$ such that, for all
$i$, $cB_i$ contains at least one point $x_i$ outside $E_\lambda$, that is
$${\cal M}(|Tf|^2)(x_i)\le \lambda^2.$$ Therefore, according to Lemma \ref{gooda} for the balls
$cB_i$, for every
 $\gamma>0$ and $K>K_0$ 
\begin{equation}
 \mu\left(U_{\lambda,i}\right)\le
C(\gamma^2 + K^{-p_0})
\,\mu(B_i),\label{goodi}
\end{equation}
where $U_{\lambda,i}=\{x\in
cB_i; {\cal M}(|Tf|^2)(x)>K^2\lambda^2,\ {\cal M}_{T,A}^\#f(x)\le
\gamma\lambda\}$.

Let 
$$U_\lambda=\{x\in
M; {\cal M}(|Tf|^2)(x)>K^2\lambda^2,\ {\cal M}_{T,A}^\#f(x)\le
\gamma\lambda\}.$$
Then, since $K>1$,
$$U_\lambda \subset E_\lambda\subset\cup_i (cB_i),$$
thus, for  all $\lambda>\lambda_0$,
\begin{eqnarray*}
\mu(U_\lambda)&=&\sum_i\mu\left(U_{\lambda,i}\right)\le C(\gamma^2 + K^{-p_0}) \sum_i\mu(B_i)\\
&\le&  C'(\gamma^2 + K^{-p_0})\,\mu(E_\lambda).
\end{eqnarray*}
Now
\begin{eqnarray*}
I_2&=&K^{p}\int_{\lambda_0}^{\infty}p\lambda^{p-1}\mu\{{\cal
M}(|Tf|^2)>K^2\lambda^2\}\,d\lambda\\
&\le&K^{p}\int_{\lambda_0}^{\infty}p\lambda^{p-1}
\left(\mu(U_\lambda)+\mu(\{{\cal M}_{T,A}^\#f>\gamma \lambda\})\right)\,d\lambda\\
\\&\le&K^{p}\int_{\lambda_0}^{\infty}p\lambda^{p-1}
\left(C'(\gamma^2 + K^{-p_0})\,\mu(E_\lambda)+\mu(\{{\cal M}_{T,A}^\#f>\gamma
\lambda\})\right)\,d\lambda\\
&=&C'(K^{p-p_0}+K^p\gamma^2)\|\left({\cal
M}(|Tf|^2)\right)^{1/2}\|_p^p+K^{p}\gamma^{-p}\|{\cal M}_{T,A}^\#f\|_p^p.
\end{eqnarray*}
Since $p<p_0$,  one obtains the lemma by choosing first $K$ large enough and then   $\gamma$ small
enough.

\bigskip

Now we are ready to prove Theorem \ref{lp}. Let $f \in L^2(M,\mu) \cap L^p(M,\mu)$. Then $Tf \in
L^p(M,\mu)$ and by Lemma \ref{sharpi}
$$
\| Tf \|_p \le  \|\left({\cal
M}(|Tf|^2)\right)^{1/2}\|_p \le C \left(\|{\cal M}_{T,A}^\#f\|_p+\|f\|_p\right).
$$
Using (\ref{pointwise}) and the strong type $(p/2,p/2)$ of the maximal function 
yields
$$
 \| Tf \|_p \le C\|f\|_p
$$
and the conclusion follows by density.

\bigskip

\noindent{\bf Remark:}  We implemented an algorithm with a right regularization of the operator $T$
by looking at 
$TA_r$ and $T(I-A_r)$. It is
also possible to obtain a result with a left regularization by making  assumptions on 
$A_r$ and $(I-A_r)T$.

If one can use duality (that is, if $T$ is linear), one can try to prove $L^p$
boundedness of
$T^*$ for some $p<2$. Then one can invoke a result in \cite{DuMcIn} if one wants a result for 
the full range $1<p<2$ or its generalization to a  limited range  $p_0<p<2$ in \cite{BK1}. 
 
Another way (which covers the sublinear case as well), would be to mimic what we just did using instead
the sharp function  introduced in \cite{M}.  Define for $f\in L^2(M,\mu)$,
$${\cal M}^{\#}_{A}f(x)=\sup_{B\ni x}\frac{1}{\mu(B)}\int_B|f-A_{r(B)}f|\,d\mu,$$
 where the supremum is taken over all balls $B$ in $M$ containing $x$,  and $r(B)$ is  the radius
of $B$. Regularizing $T$ from the left means considering ${\cal M}^{\#}_{A}(Tf)$. 

However, if we were to apply left regularization to the boundedness of Riesz transforms on
manifolds, since $T$ takes functions to vector fields (or to $1$-forms), $A_r$ would have to be
vector-valued, say, the heat semigroup on $1$-forms instead of the heat semigroup on functions
as in Section \ref{Riesz} below. Since only the action of $A_r$ on the image  of $T$ comes into play, one sees that 
the needed assumptions on  the heat semigroup on $1$-forms would only concern its action on {\em exact} forms;
this fits with our purpose, see the discussion at the end  of Section \ref{assumption}.
The advantage of the right regularization is that it is susceptible to be applied to more general situations, where the notion
of differential  forms is not available. 

 For more on the applications of the above sharp function, the associated space, and its
potential uses in harmonic analysis, see \cite{DY}.

\subsection{The local criterion}\label{local}

The theorem above admits variations towards localization. 

Denote by  ${\cal M}_E$  the 
Hardy-Littlewood maximal operator relative to a measurable subset $E$ of $M$, that is, for $x\in E$ and $f$ a locally
integrable function on $M$,
$${\cal M}_Ef(x)=\sup_{B\ni x}\frac{1}{\mu(B \cap E)}\int_{B\cap E}|f|\,d\mu, $$
where $B$ ranges over all open balls of $M$ containing $x$ and centered in $E$. If in particular, 
$E$ is a ball with radius $r$, it is enough to consider  balls $B$ with radii not exceeding $2r$.

We say that a subset $E$ of $M$ has the relative doubling property if there exists a constant $C_E$ such
that  for all $x \in E$ and $r>0$ we have 
$$
\mu(B(x,2r)\cap E)  \le C_E \mu(B(x,r)\cap E).
$$
In other words, $E$ endowed with the induced distance and measure has the doubling property. 
The constant $C_E$ is called the relative doubling constant of $E$. On such a set, ${\cal M}_E$ is 
weak type $(1,1)$ and bounded on $L^p(E,\mu)$, $1<p\le \infty$.

\begin{theo}\label{lploc} Let $(M, d,\mu)$ be a measured metric space. Let $p_0\in (2,\infty]$.
Suppose
that $T$ is a bounded sublinear operator which is bounded on $L^2(M,\mu)$, and let  $A_r$, $r>0$,
be  a family of linear operators acting on $L^2(M,\mu)$.
Let $E_1$ and $E_2$ be two subsets of $M$ such that $E_2$ has the relative doubling property,
$\mu(E_2)<\infty$ and
$E_1\subset E_2$. 
     Assume  that $f \mapsto {\cal M}_{E_2,T,A}^\#f$ is bounded from $L^p(E_1, \mu)$ into
$L^p(E_2,\mu)$   for all $p \in (2, p_0)$, where
\begin{equation}
\big({\cal M}_{E_2,T,A}^\#f\big)^2(x) = \sup_{B\ {ball\ in }\ M, \ B \ni x\ }\frac{1}{\mu(B\cap
E_2)}\int_{B\cap E_2}|T(I-A_{r(B)})f|^2\,d\mu, \quad x \in E_2,\label{pointwiseloc}
\end{equation}
and, for some  sublinear operator $S$ bounded from $L^p(E_1, \mu)$ into
$L^p(E_2,\mu)$   for all $p \in (2, p_0)$, 
\begin{equation}
\left(\frac{1}{\mu(B\cap E_2)}\int_{B\cap E_2} |TA_{r(B)}f|^{p_0} \,d\mu\right)^{1/p_0}\le 
C\big({\cal M}_{E_2} (|Tf|^2)+ (Sf)^2\big)^{1/2} (x),\label{domaloc}
\end{equation}
for all $f\in L^2(M,\mu)$ supported in $E_1$, all balls $B$ in $M$ and all $x \in B\cap E_2$, where $r(B)$
is the radius of
$B$. If  $2<p<p_0$ and  $Tf \in L^p(E_2,\mu)$ whenever $f\in L^p(E_1,\mu)$,  then $T$ is bounded from
$L^p(E_1,\mu)$  into $L^p(E_2,\mu)$ and its operator norm is bounded by a constant depending only on
the operator norm of $T$ on $L^2(M,\mu)$,  $C_{E_2}$,
$p$, $p_0$, the operator norms of ${\cal M}_{E_2,T,A}^\#$ and $S$ on $L^p$, and the constant in \eqref{domaloc}. 
\end{theo}

Again, if $p_0=\infty$ the left-hand side of \eqref{domaloc} is understood as the essential
supremum on $B$. 

The proof of this result is almost identical to that of Theorem \ref{loc} once we make some
adjustments. The first one is to forget about $M$ and to work directly in the relative space $E_2$
by replacing systematically
$T$ and
$TA_r$ by truncations
$\chi_{E_2}T\chi_{E_1}$ and 
$\chi_{E_2}TA_r\chi_{E_1}$. Thus the maximal operator relative to  $E_2$ becomes the maximal operator
on $E_2$. 

The second one is that Lemma \ref{gooda}  becomes

\begin{lem}\label{goodaloc} Let $(M, d,\mu, A_r, p_0, E_2, E_1)$ and $T$ be as above. 
Assume  that 
${\rm (\ref{domaloc})}$ holds. There exist $K_0>1$ and $C>0$ only depending on $C_{E_2}$ and $p_0$,
such that, for every $\lambda>0$, every
$K>K_0$ and $\gamma>0$, for every ball $B_0$ in $M$ and every function $f\in L^2(E_1,\mu)$ such
that there exists $x_0\in B_0\cap E_2$ with ${\cal M}_{E_2}
(|Tf|^2)(x_0)+ (Sf)^2 (x_0) \le \lambda^2$, then 
$$
\mu\left(\{x\in B_0\cap E_2; {\cal M}_{E_2}(|Tf|^2)(x)>K^2\lambda^2,\ {\cal M}_{E_2,T,A}^\#f(x)\le
\gamma\lambda\}\right)\le
C(\gamma^2 + K^{-p_0})
\,\mu(B_0\cap E_2).
$$
\end{lem}
 The proof is the same
since $Sf\ge 0$ implies ${\cal M}_{E_2}
(|Tf|^2)(x_0) \le \lambda^2$.

The third one is that the term $Sf$ brings a modification in Lemma \ref{sharpi} which becomes

\begin{lem}\label{sharpiloc} Let $(M, d,\mu, A_r, E_2, E_1)$  and $T$ be as above. 
Assume  that 
 ${\rm (\ref{domaloc})}$ holds.
Then, for  $0< p< p_0$, 
$$
\|\left({\cal M}_{E_2}(|Tf|^2)\right)^{1/2}\|_{L^p(E_2)}\le C\big(\|{\cal
M}_{E_2,T,A}^\#f\|_{L^p(E_2)}+
\|Sf\|_{L^p(E_2)}+ \|f\|_{L^p(E_2)}\big),
$$
for every $f\in
L^2(E_1,\mu)$ for which the left-hand side is finite, where $C$ depends only on $p, p_0$ and the
doubling constant of $E_2$ (but not on $\mu(E_2)$).  
\end{lem} 

Since $\mu(E_2)<\infty$, the Whitney decomposition in the proof of Lemma \ref{sharpi} can be
performed  for
$\lambda>\lambda_0$ with $\lambda_0=
\frac{1}{\mu(E_2)} \int_{E_2}  {\cal M}_{E_2}(|Tf|^2)^{1/2}\, d\mu$. 
Again, by Kolmogorov inequality, the doubling property of $E_2$,  and
the weak type
$(1,1)$ of
${\cal M}_{E_2}$ (with constant independent of the size of $\mu(E_2)$),
$$
\int_{E_2}  {\cal M}_{E_2}(|Tf|^2)^{1/2}\, d\mu \le C\mu(E_2)^{1/2} \||Tf|^2 \|_{L^1(E_2)}^{1/2} =
C\mu(E_2)^{1/2} \|Tf \|_{L^2(E_2)}.
$$
Using then the $L^2$ boundedness of $T$ and the support condition of $f$ (this is where we use
$E_1\subset E_2$), we obtain
$$
\lambda_0 \le \frac{C}{\mu(E_2)^{1/2}} \|f\|_{L^2(E_2)} \le  \frac{C}{\mu(E_2)^{1/p}}
\|f\|_{L^p(E_2)}.
$$
Now write with the notation of the proof in Lemma \ref{sharpi}  
\begin{align*}\|\left({\cal M}_{E_2}(|Tf|^2)\right)^{1/2}\|_{L^p(E_2)}^p& = \int_{{\cal
M}_{E_2}(|Tf|^2)
\le K^2\lambda_0^2} \left({\cal M}_{E_2}(|Tf|^2)\right)^{p/2}\, d\mu \\
& \qquad \qquad + 
\int_{{\cal M}_{E_2}(|Tf|^2)
>K^2\lambda_0^2} \left({\cal M}_{E_2}(|Tf|^2)\right)^{p/2} \, d\mu.
\end{align*}
The last integral can be treated as before, using Lemma \ref{goodaloc}. The first integral 
is bounded above by 
$$
K^p\lambda_0^p \mu(E_2) \le K^p C^p \|f\|_{L^p(E_2)}^p
$$
with $C$ depending only on the $L^2$ norm of $T$ and the doubling constant of $E_2$. Further details are left to the reader.

\section{Application to the Riesz transform}\setcounter{equation}{0}\label{Riesz}

  In this section, $M$ is a complete non-compact Riemannian manifold,
 $\Delta$ denotes the Laplace-Beltrami operator,
$e^{-t\Delta}$, $t>0$, the heat semigroup and $p_t(x,y)$, $t>0$, $x,y\in M$, the heat kernel. The measure $\mu$ is the induced 
Riemannian volume. The measure of the ball $B(x,r)$, $x\in M$, $r>0$ is
also written $V(x,r)$. 

We prove the statements corresponding to the (global) Riesz transform $\nabla  \Delta^{-1/2}$.
We set $Tf=|\nabla  \Delta^{-1/2}f|$ (remember that, in the finite volume case, we restrict ourselves to
functions with mean zero; in other words, in order to apply verbatim Theorem
\ref{lp}, we define $\nabla  \Delta^{-1/2}$ by zero on constants). The boundedness of $T$ on $L^2$ has been already observed. 

\subsection{Proof of the main result}

We
now show    Theorem
\ref{mainp>2},  namely the fact that, under 
$(D)$ and
$(P)$,   $(G_{p_0})$ implies the $L^p$ boundedness of the  Riesz transform for $2<p<p_0$, as a  consequence of Theorem
\ref{lp}. 

Recall that Theorem
\ref{lp} applies if $T$ is assumed to act  on $L^p(M)$. However, we have 
$$
\nabla \Delta^{-1/2} = c\int_0^\infty  \nabla e^{-t\Delta} \,\frac{dt}{\sqrt t}.
$$
If we set  $T_\ep f=| c\int_\ep^{1/\ep} \nabla e^{-t\Delta}f \,\frac{dt}{\sqrt t}|
$ for $0<\ep<1$, then for $f\in L^2(M)$ we have  $\|T_\ep f\|_2 \le \|f\|_2$ (this follows from (\ref{triv})
and spectral theory) and    $T_\ep f
\to Tf$ in
$L^2(M)$ as
$\ep
\to 0$, while $\|T_\ep f\|_p \le C_\ep \|f\|_p$ for $f \in  L^p(M)$. As the application of Theorem
\ref{lp} to $T_\ep$ gives us a uniform bound with respect to $\ep$, a limiting argument yields 
the $L^p$ boundedness of $T$ on $L^2(M) \cap L^p(M)$, hence on $L^p(M)$. Henceforth, we ignore this
approximation step and our goal is now to establish  (\ref{pointwise}) and
(\ref{doma}) for $T$. 
\bigskip

The first ingredient is Gaffney off-diagonal estimates valid in a general Riemannian manifold:
There exist two constants
$C\ge 0$ and $\alpha>0$  such that,  for every
$t>0$,  every closed subsets $E$ and $F$ of $M$, and every function $f$ supported in $E$, one has
\begin{equation}\label{Gaffney}
\| e^{-t\Delta} f \|_{L^2(F)} + \| t\Delta e^{-t\Delta} f \|_{L^2(F)} +
\|\sqrt t \, |\nabla e^{-t\Delta} f| \|_{L^2(F)}  
 \le C e^{-\frac {\alpha d(E,F)^2} t} \|f\|_{L^2(E)}.
\end{equation}
Here, $d(E,F)$ is the distance between  the sets $E$ and $F$. 
The inequality for the first  term of the left-hand side is classical (see, e.g., \cite{Dav}).  The estimate for the second
one follows from essentially the same proof (see \cite{Davhigher}, Lemma 7).  We give a proof for the third one as we could not
find it in the literature in this situation (it is proved in
\cite{AT} for elliptic operators on $\mathbb{R}^n$).

We assume $d(E,F) > \sqrt t$ as otherwise there is nothing to prove. Let  $\tilde F$ be
the set of those
$x\in M$ for which  $d(x,F)\leq\frac {d(E,F)}{3}$.  Let $\varphi$ be a smooth function on $M$ such that $0 \le \varphi \le
1$,  
$\varphi$ is supported in 
$\tilde F$, $\varphi\equiv 1$ on
$F$, and $|\nabla \varphi| \le \frac 6 {d(E,F)}$.  Set
$A=\|\sqrt t \, \varphi |\nabla e^{-t\Delta} f|
\|_{2}\ge \|\sqrt t
\, |\nabla e^{-t\Delta} f|
\|_{L^2(F)}$.   Integrating by parts, 
\begin{align*}
A^2 &= t \langle \varphi^2 \nabla e^{-t\Delta} f, \nabla e^{-t\Delta} f \rangle \\
&  = - 2t \langle   (e^{-t\Delta} f)\nabla \varphi, \varphi \nabla  e^{-t\Delta} f \rangle 
+ t \langle \varphi^2  e^{-t\Delta} f, \Delta e^{-t\Delta} f \rangle \\
& \le 2\sqrt{t} \| |\nabla \varphi| e^{-t\Delta} f \|_2 A +  \|\varphi e^{-t\Delta} f \|_{2} \|
t \, \varphi\Delta e^{-t\Delta} f \|_{2}.
\end{align*}
Using the  properties of $\varphi$ and the bounds for the first two terms in \eqref{Gaffney} 
we obtain the desired conclusion. 

\bigskip

We now introduce the regularizing operator $A_r$, $r>0$,  by setting $$I-A_r= (I-e^{-r^2\Delta})^n$$
for some integer $n$ to be chosen. Observe that $A_r$ is bounded on $L^2(M)$
with norm 1. We prove  (\ref{pointwise}) in the following lemma.  

\begin{lem}\label{lempointwise} Assume that $(D)$ holds. Then, for some $n$ large enough depending
only on
$(D)$, for every ball $B$ with radius $r>0$ and all $x\in B$,  
\begin{equation}\label{pointwiseriesz}
\| |\nabla \Delta^{-1/2}
(I-e^{-r^2\Delta})^nf|\|_{L^2(B)} \le C  \mu(B)^{1/2}  \big({\cal
M}(|f|^2)(x)\big)^{1/2}.
\end{equation}
  
\end{lem}

\noindent{\bf Proof:} Let $f \in L^2(M)$. Take a ball $B$ with radius $r=r(B)$  and $x$ a point in
$B$.  Denote by
$C_i$ the ring $2^{i+1}B\setminus 2^iB$ if $i\ge 2$ and let $C_1=4B$. Decompose $f$ as $f_1+f_2+f_3+\ldots$
with 
$f_i=f\chi_{C_i}$. By Minkowski inequality we have that 
$$
\| |\nabla \Delta^{-1/2} (I-e^{-r^2\Delta})^nf|\|_{L^2(B)} \le \sum_{i\ge 1} 
\| |\nabla \Delta^{-1/2}
(I-e^{-r^2\Delta})^nf_i|\|_{L^2(B)}.
$$
For $i=1$ we use the $L^2$ boundedness of $\nabla \Delta^{-1/2}
(I-e^{-r^2\Delta})^n$:
$$
\| |\nabla \Delta^{-1/2}
(I-e^{-r^2\Delta})^nf_1|\|_{L^2(B)} \le \| f \|_{L^2(4B)} \le \mu(4B)^{1/2}
 \big({\cal M}(|f|^2)(x)\big)^{1/2}.
$$
For $i\ge 2$ we use the integral representation of $\Delta^{-1/2}$:
\begin{eqnarray*}\nabla \Delta^{-1/2}
(I-e^{-r^2\Delta})^n&=&c\int_0^\infty  \nabla e^{-t\Delta} (I-e^{-r^2\Delta})^n\,\frac{dt}{\sqrt
t}\\ &=&c\int_0^\infty  g_{r}(t) \nabla e^{-t\Delta} \, {dt}\\
\end{eqnarray*}
where using the usual notation for the binomial coefficient, 
$$g_r(t)=\sum_{k=0}^n \binom n k (-1)^k \frac{\chi_{\{t>kr^2\}}}
{\sqrt{t-kr^2}}.$$  

By Minkowski integral inequality and Gaffney estimates (\ref{Gaffney}), using the support of $f_i$, we
have that 
$$
\| |\nabla \Delta^{-1/2}
(I-e^{-r^2\Delta})^nf_i|\|_{L^2(B)} \le C \left(\int_0^\infty |g_{r}(t)| e^{-\frac
{\alpha' 4^ir^2} t}  \frac {dt} {\sqrt t} \right) \| f\|_{L^2(C_i)}
$$ 
The latter integral can be estimated as follows.  
Elementary analysis yields the following estimates for $g_{r}$:
$$
|g_{r}(t)| \le \frac {C_n} {\sqrt{t-\ell r^2}} \quad \textrm{if} \quad {0\le \ell r^2<t \le (\ell+1)r^2 \le (n+1)r^2}
$$
and 
$$
|g_{r}(t)| \le C_nr^{2n}t^{-n-\frac{1}{2}} \quad \textrm{if} \quad{t> (n+1)r^2}.
$$
The latter estimate comes from the inequality 
$$\left|\sum_{k=0}^n \binom n k (-1)^k \varphi(t-kr^2) \right| \le C_n\sup_{u\ge \frac{t}{n+1}} |\varphi^{(n)}(u)| r^{2n},
$$
 which can be obtained by expanding $\varphi(t-ks)$ using Taylor's formula about $t$ and using the
classical  relations $\sum_{k=0}^n \binom n k (-1)^k k^\ell =0$ for $\ell=0,...,n-1$ (see \cite{F}, problem 16, p.65).  
This yields the following estimates, uniformly  with respect to $r$: 
$$
\int_0^\infty |g_{r}(t)| e^{-\frac
{\alpha' 4^ir^2} t}  \frac {dt} {\sqrt t} \le C_n 4^{-in}.
$$
Now, an easy consequence of $(D)$ is that  for all $ y\in M,  
r>0,$ and $ \theta\ge 1$ 
\begin{equation}
V(y,\theta r)\le  C\theta^\nu V(y,r),\label{nu}
\end{equation}
 for some constants $C$ and $\nu>0$.
Therefore, since $C_i \subset 2^{i+1}B$, 
$$\| f\|_{L^2(C_i)}  \le \mu(2^{i+1}B)^{1/2}
 ({\cal M}(|f|^2)(x))^{1/2} \le  \sqrt{C}2^{(i+1)\nu/2} \mu(B)^{1/2}  \big({\cal M}(|f|^2)(x)\big)^{1/2}.
$$
Choosing $2n>\nu/2$, we have 
$$
\| |\nabla \Delta^{-1/2}
(I-e^{-r^2\Delta})^nf|\|_{L^2(B)} \le C' \left(\sum_{i\ge 1}  2^{i(\nu/2-2n)}\right) \mu(B)^{1/2}  \big({\cal
M}(|f|^2)(x)\big)^{1/2},
$$
which proves Lemma \ref{lempointwise}. 
\bigskip

We now show that (\ref{doma}) holds. We begin with a  lemma. 

\begin{lem}\label{lemoffdiagl2lp} Assume  $(D)$, $(P)$ and $(G_{p_0})$. Then  the following estimates hold: for
every
$p\in (2,p_0)$, for every ball
$B$ with radius
$r$ and every
$L^2$ function $f$ supported in $C_i=2^{i+1}B\setminus 2^iB$, $i\ge 2$, or $C_1=4B$, and every $k
\in
\{1, 
\ldots, n\}$, where $n$ is chosen as above, one has
\begin{equation}
 \left( \frac 1{\mu(B)}\int_B   |\nabla e^{-kr^2\Delta}f|^p\, d\mu \right)^{1/p}   \le 
\frac {Ce^{- \alpha 4^i}}{r}\left(\frac 1 {\mu(2^{i+1}B)} \int_{C_i} |f|^2\, d\mu\right)^{1/2}
\label{symi}
\end{equation}
for some constants $C$ and $\alpha$ depending only on $(D)$, $(P)$,  $p$ and $p_0$.
\end{lem}

\noindent{\bf Proof:} 
 By interpolating $(G_{p_0})$ with 
the $L^2$ Gaffney estimates, we obtain  $L^p$ Gaffney estimates  for any $p\in (2,p_0)$:
for every
$t>0$, for every closed sets $E$ and $F$ and every function $f$ supported in $E$, one has
\begin{equation}\label{Gaffneyp}
 \|\sqrt t \, |\nabla e^{-t\Delta} f| \|_{L^p(F)}  \le C e^{-\frac {\alpha d(E,F)^2}
t}\|f\|_{L^p(E)},
\end{equation} 
with $C>0$ and $\alpha>0$ depending on $p,p_0$, and the constant  $C_{p_0}$ in $(G_{p_0})$.
Now let   $B$ be a ball with  radius  $r$ and let $f$ be 
supported in $C_i$.

 In the following proof, many constants will implicitely depend  on $n$, which itself only depends on
$(D)$.

Let us begin with the case $i=1$. The above estimate (or, directly,
$(G_p)$) yields 
\begin{equation}
\left(\int_B   |\nabla e^{-kr^2\Delta}f|^p\, d\mu\right)^{1/p} \le \frac {C}{r} \left(\int_M
|e^{-(k/2)r^2\Delta}f|^p \, d\mu\right)^{1/p}.
\label{sym1grad}
\end{equation}
Let $t=(k/2)r^2$. Since $(U\!E)$ follows from  $(D)$ and $(P)$,  one has the upper
estimate
$$p_t(x,y)
\le
\frac{C} {V(y,\sqrt{t}) }\exp\left(-c\frac{d^2(x,y)}{t}\right),$$
for all $x,y\in M$.
Because of the doubling property, $$ V(y,\sqrt{t}) \simeq V(x_B,\sqrt{t})\simeq \mu(B)$$ where $x_B$ is the
center of
$B$ and
$y\in 4B$. It follows that
\begin{equation}
| e^{-(k/2)r^2\Delta}f(x)|\le  
\left(\frac C{\mu(B)}\int_{4B} |f|^2\, d\mu\right)^{1/2},
\label{sym1linfty}
\end{equation}
for all $x\in M$.
On the other hand, by Gaffney (or $L^2$ contractivity of the heat semigroup),
\begin{equation}
\int_{M} | e^{-(k/2)r^2\Delta}f|^2\, d\mu   \le  C\int_{4B} |f|^2\, d\mu.\label{sym12}
\end{equation}
Thus, by H\"older,  
$$ \left(\int_M
|e^{-(k/2)r^2\Delta}f|^p\, d\mu \right)^{1/p} \le C \mu(B)^{\frac 1 p - \frac 1 2} \left(\int_{4B}
|f|^2\, d\mu\right)^{1/2},
$$
which, together with \eqref{sym1grad}, yields \eqref{symi}  in this case. 

Next assume that $i\ge 2$. 
Denote by $\chi_{C\ell}$ the characteristic function of $C_\ell$ and write
$$
\nabla e^{-kr^2\Delta} f = \sum_{\ell\ge 1} h_\ell, \quad h_\ell=\nabla e^{-(k/2)r^2\Delta}
(\chi_{C\ell}) e^{-(k/2)r^2\Delta}f.
$$
 From (\ref{Gaffneyp})  we have
$$
\left( \frac 1{\mu(B)}\int_B   |h_\ell|^p\, d\mu
\right)^{1/p}
\le \left( \frac{\mu(2^{\ell+1}B)}{\mu(B)}\right)^{1/p} \frac {Ce^{- \alpha 4^\ell}}{r} 
\left(\frac 1{\mu(2^{\ell+1}B)}\int_{C_\ell} | e^{-(k/2)r^2\Delta}f|^p\, d\mu\right)^{1/p}
$$
and, using  (\ref{nu}),
\begin{equation}
\left( \frac 1{\mu(B)}\int_B   |h_\ell|^p\, d\mu
\right)^{1/p} \le C'2^{(\ell+1)\nu/p} \frac {e^{- \alpha 4^\ell}}{r} 
\left(\frac 1{\mu(2^{\ell+1}B)}\int_{C_\ell} | e^{-(k/2)r^2\Delta}f|^p\, d\mu\right)^{1/p}.\label{symgrad}
\end{equation}

From the $L^2$ Gaffney
estimates for the semigroup, one has 
$$\int_{C_\ell} | e^{-(k/2)r^2\Delta}f|^2\, d\mu   \le  K_{i\ell} \int_{C_i} |f|^2\, d\mu
$$
with  
$$
K_{i\ell} \le   \begin{cases} C e^{-c4^i} &\text{if} \ \ell\le i-2,\\
C  &\text{if} \ i-1\le \ell \le i+1,\\
Ce^{-c4^\ell}  & \text{if} \ \ell \ge i+2. \\
\end{cases}
$$
Since, by (\ref{nu}), $K_{i\ell} \frac{\mu(2^{i+1}B)}{\mu(2^{\ell+1}B)} \le K_{i\ell} \sup(1, C2^{(i-\ell)\nu}) $, and if we
still denote by $K_{i\ell}$  a sequence of the same  form  with different constants,  we may also write
\begin{equation}
\frac 1{\mu(2^{\ell+1}B)}\int_{C_\ell} | e^{-(k/2)r^2\Delta}f|^2\, d\mu   \le  K_{i\ell}\ \frac
1{\mu(2^{i+1}B)}\int_{C_i} |f|^2\, d\mu.\label{sym2}
\end{equation}
Next, it easily follows from $(U\!E)$ and $(D)$ that for all  $x \in C_\ell$, 
\begin{align*} 
| e^{-(k/2)r^2\Delta}f(x)|  &\le  C\int_{C_i}{V(y,r)}^{-1}
\exp\left(-c\frac{d^2(x,y)}{r^2}\right) |f(y)| \, d\mu(y)
\\
& \le  K_{i\ell} \int_{C_i}{V(y,r)}^{-1} |f(y)| \, d\mu(y).
\end{align*}
If $y \in C_i$,   $2^{i+1} B \subset B(y, 2^{i+2}r)$, so that 
$$  \frac 1 {{V(y,r)}} \le  \frac {C2^{(i+2)\nu}}{ {V(y,2^{i+2}r)}} \le
\frac {C 2^{(i+2)\nu}}{ \mu(2^{i+1}B)},$$
and it follows that 
\begin{equation}
| e^{-(k/2)r^2\Delta}f(x)|\le  K_{i\ell}  2^{(i+2)\nu} 
\left(\frac 1{\mu(2^{i+1}B)}\int_{C_i} |f|^2\, d\mu\right)^{1/2}
\label{symlinfty}
\end{equation}

By applying H\"older and using  \eqref{sym2} and \eqref{symlinfty}, one obtains
$$
\left(\frac 1{\mu(2^{\ell+1}B)}\int_{C_\ell} | e^{-(k/2)r^2\Delta}f|^p\, d\mu\right)^{1/p}   \le 
K_{i\ell} 2^{(i+2)\nu (1 -\frac{2}{p})} \left(
\frac 1{\mu(2^{i+1}B)}\int_{C_i} |f|^2\, d\mu\right)^{1/2}.$$
Together with \eqref{symgrad} and summing in $\ell$, this yields  
\begin{align*}
\left( \frac 1{\mu(B)}\int_B   |\nabla e^{-kr^2\Delta}f|^p\, d\mu \right)^{1/p} &\le C \sum_{\ell \ge 1}
2^{(\ell+1)\nu/p}
\frac{e^{-c4^\ell}}{r} K_{i\ell}
2^{(i+2)\nu (1 -\frac{2}{p})}
\left(\frac 1{\mu(2^{i+1}B)}\int_{C_i} |f|^2\, d\mu\right)^{1/2}
\\
& \le C \frac{e^{-c'4^i}}{r} \left(\frac 1{\mu(2^{i+1}B)}\int_{C_i} |f|^2\, d\mu\right)^{1/2}.
\end{align*}
  This ends the proof of Lemma
\ref{lemoffdiagl2lp}. 
\bigskip

Equipped with this lemma, we can prove (\ref{doma}) for any $p\in (2,p_0)$. 
Fix such a $p$.  By expanding $I-(I-e^{-r^2\Delta})^n$ it suffices to show
\begin{equation}\label{maxp}
\left(\frac{1}{\mu(B)}\int_B |\nabla e^{-kr^2\Delta}f|^{p} \,d\mu\right)^{1/p}\le  C\big({\cal M}
(|\nabla f|^2)\big)^{1/2} (y)
\end{equation}
for  $f $ with $f, \nabla f$ locally square integrable, $B$ any ball with $r=r(B)$, $y\in B$ and 
$k=1,2, \ldots, n$. Recall that $n$ is  chosen larger than $\nu/4$ where $\nu$ is given in
(\ref{nu}).

Recall that our assumptions ensure  that $M$ satisfies (\ref{stoco}). In other words,
$$e^{-t\Delta}1\equiv  1,\ \forall\,t>0.$$   

We may therefore write
$$
\nabla e^{-kr^2\Delta}f = \nabla e^{-kr^2\Delta}(f- f_{4B}).$$
Write $f- f_{4B}= f_1+f_2+f_3+\ldots$ where 
 $f_i= (f- f_{4B})\chi_{C_i}$.
For $i=1$, we use the lemma and $(P)$ to obtain
$$
\left( \frac 1{\mu(B)}\int_B   |\nabla e^{-kr^2\Delta}f_1|^p\, d\mu \right)^{1/p} 
\le C\left(\frac {1} {\mu(4B)} \int_{4B}
|\nabla f|^2\,d\mu \right)^{1/2} \le C  \big({\cal M}
(|\nabla f|^2)\big)^{1/2} (y).
$$
For $i\ge 2$, we have similarly
$$
\left( \frac 1{\mu(B)}\int_B   |\nabla e^{-kr^2\Delta}f_i|^p\, d\mu \right)^{1/p} 
\le \frac {Ce^{- \alpha 4^i}}{r} \left(\frac {1} {\mu(2^{i+1}B)} \int_{C_i}
|f_i|^2\,d\mu \right)^{1/2}.
$$
But 
$$
\int_{C_i}
|f_i|^2\,d\mu \le \int_{2^{i+1}B}
|f - f_{4B}|^2\,d\mu,
$$
$$
|f - f_{4B}| \le |f-f_{2^{i+1}B}|
+\sum_{\ell=2}^{i}|f_{2^\ell B}-f_{2^{\ell+1}B}|
$$
and observe that 
$$
|f_{2^\ell B}-f_{2^{\ell+1}B}|^2  \le \frac 1 {\mu(2^{\ell+1}B)} \int_{2^{\ell+1}B}
|f - f_{2^{\ell+1}B}|^2\,d\mu \le   (2^{\ell}r)^2 {\cal M}
(|\nabla f|^2) (y).
$$
Hence, by Minkowski inequality, we easily obtain
\begin{equation}\label{eqci}
\left(\frac {1} {\mu(2^{i+1}B)} \int_{C_i}
|f_i|^2\,d\mu \right)^{1/2} \le  C  (2^{i}r)
\big({\cal M}
(|\nabla f|^2)\big)^{1/2} (y).
\end{equation}
It remains to sum for $i\ge 2$
\begin{align*}
\sum_{i\ge 2}\frac {Ce^{- \alpha 4^i}}{r} \left(\frac {1} {\mu(2^{i+1}B)} \int_{C_i}
|f_i|^2\,d\mu \right)^{1/2} &\le 
\sum_{i\ge 2} \frac {Ce^{- \alpha 4^i}}{r}  ( 2^{i}r ) \big({\cal M}
(|\nabla f|^2)\big)^{1/2} (y)\\
&\le C \sum_{i\ge 2} {e^{- \alpha 4^i}} 2^{i} \big({\cal M}
(|\nabla f|^2)\big)^{1/2} (y).
\end{align*}
This yields \eqref{maxp}, and the proof of Theorem \ref{mainp>2} is finished.

\bigskip

\noindent{\bf Remark:} 

-A slight modification of the proof allows the following improvement of 
\eqref{maxp}: for some constants $c,C>0$
\begin{equation}\label{maxpimp}
\left(\frac{1}{\mu(B)}\int_B |\nabla e^{-kr^2\Delta}f|^{p} \,d\mu\right)^{1/p} \le  C \inf_{x
\in B}\big(  e^{-cr^2\Delta} (|\nabla f|^2)\big)^{1/2} (x)
\end{equation}
for  $f $ with $f, \nabla f$ locally square integrable, $B$ any ball with $r=r(B)$ and 
$k=1,2, \ldots, n$, where $n$ is chosen as above. See Lemma \ref{Steve} in Section \ref{sectionsteve} for    the case
$p=\infty$ of this inequality.

-The techniques of proof of Theorem \ref{lp} extend easily to the
vector-valued setting. It can be checked that it applies as well to obtain the square function
estimate (seen as
$L^p$ boundedness of a vector-valued operator)
\begin{equation} \label{squarefunction}
\left\| \left(\int_0^\infty |\nabla e^{-t\Delta}f|^2\, {dt}\right)^{1/2} \right\|_p \le C
\|f\|_p
\end{equation}
from $(D)$, 
$(P)$, and $(G_{p_0})$, for $2<p<p_0$.  Compare with \cite{CDfull}, Section 3, where it is observed  that a related inequality
(with the Poisson semigroup instead of the heat semigroup) holds under
\eqref{conj} only. In particular, $(D)$ and $(P)$ are not used there.  This explains further our
discussion in Section
\ref{assumption} where we question the relevance of $(D)$ and $(P)$ for the Riesz transform
boundedness. 

\subsection{A simpler situation}\label{sectionsteve}

In this section, we show a slightly weaker version of Theorem \ref{maincor}  by replacing
$(G)$ with the stronger inequality
\eqref{gradgub}, that is
\begin{equation*} 
|\nabla_x\, p_t(x,y)|\le \frac{C}{\sqrt{t}\,V(y,\sqrt{t})}\exp{\left(-c\frac{d^2(x,y)}{t}\right)},\ \forall\,x,y\in M,\ t>0,
\end{equation*}
as this is a simpler application of Theorem \ref{lp} and this 
 hypothesis is related to the domination condition \eqref{conj}, therefore to the conjecture in \cite{CDfull} explained in
Section
\ref{assumption}.

We set again $Tf=|\nabla  \Delta^{-1/2}f|$ but choose here 
$A_r=e^{-r^2\Delta}$, and apply Theorem \ref{lp} in the
case $p_0=\infty$. 

\bigskip
We begin with the verification of \eqref{pointwise}.  
Let $B$ be a ball of radius $r$. Let $f \in L^2(M)$. Write $f=f_1+f_2$  where $f_1=f$ on $2B$ and 0 elsewhere.
First, the $L^2$ boundedness of $T(I-A_r)$ gives us
$$
\int_B |T(I-A_r)f_1|^2 \, d\mu \le \int_{M} |f_1|^2\, d\mu  \le \mu(2B) {\cal M}(|f|^2)(z)
$$
whenever $z\in B$. To conclude the proof of  \eqref{pointwise}, it remains to obtain the same bound for $f_2$.
One has
\begin{eqnarray*}T(I-A_r)f_2(x)&=&|\nabla 
\Delta^{-1/2}(I-e^{-r^2\Delta})f_2(x)|\\
&=&|\int_M{\tilde k}_{r}(x,y)f_2(y)\,d\mu(y)|\\
&\le&
\int_M|{\tilde k}_{r}(x,y)||f_2(y)|\,d\mu(y),
\end{eqnarray*}
 where  
$${\tilde k}_r(x,y)=\int_0^{\infty}g_r(s)\nabla_x\,p_s(x,y)\,ds$$
and $$g_r(t)=\frac {1}{ \sqrt{t}}-\frac{\chi_{\{t>r^2\}}}{ \sqrt{t-r^2}}.$$
Using \eqref{gradgub} and following the case $n=1$
in the proof of Lemma \ref{lempointwise},
one  obtains that  if $d(x,y) \ge r$ then 
$$
|\tilde k_r(x,y) | \le \frac{C}{V(x,r)} 
\exp{\left(-c\frac{d^2(x,y)}{r^2}\right)}.
$$
This estimate and the support property of $f_2$ ensure 
$$
\int_M|{\tilde k}_{r}(x,y)||f_2(y)|\,d\mu(y) \le C {\cal M}(|f_2|)(z) \le  C {\cal M}(|f|)(z)
$$
(see for instance \cite{DR}, Proposition 2.4), therefore the pointwise
bound
$$|T(I-A_r)f_2(x)| \le  C {\cal M}(|f|)(z)$$ 
whenever $x$ and $z$ belong to $B$ , hence the bound in $L^2$ average. 

\bigskip

The next step is the verification of (\ref{doma}), which in this case becomes the maximal estimate 
\begin{equation}
\sup_{y\in B}|\nabla e^{-r^2(B)\Delta}f(y)|^2\le C' \inf_{x\in B}{\cal M} (|\nabla f|^2)(x),
\label{domna}\end{equation}
for all  balls $B$ and functions $f$ with $f, \nabla f$ square integrable. 
But, under $(D)$ and $(DU\!E)$, this follows  
from \eqref{gradgub}.  

Indeed, assume \eqref{gradgub}.  Recall from  the discussion in Sections \ref{results} and \ref{assumption} that in such
a case, the full
$(LY)$ estimates hold and that \eqref{gradgub} is equivalent to \eqref{gradptpt}, that is 
\begin{equation*} 
|\nabla_xp_t(x,y)|\le \frac{C}{ \sqrt{t}}\,p_{C't}(x,y)
\end{equation*} 
for some constants $C,C'>0$.  Also, by Cauchy-Schwarz, $(P)$ implies
\begin{equation}
\int_{B}|f-f_{B}|\,d\mu\le Cr(B)\sqrt{\mu(B)}\left(\int_{B}
|\nabla f|^2\,d\mu
\right)^{1/2}\label{P2}
\end{equation}
for any ball $B$ in $(M,d)$ with radius $r(B)$ and any $f$  with $f,\nabla f$ square integrable on 
$B$  (in fact $(P)$ and (\ref{P2}) are equivalent, see \cite{HaK0}).

Fix $B$ a ball in $M$,  $x,y\in B$, and set $t=r^2(B)$. Using \eqref{stoco} again, write
\begin{equation}
|\nabla e^{-t\Delta}f(y)|=|\nabla e^{-t\Delta}(f-f(x))(y)|
\le\int_M|\nabla p_t(y,z)||f(z)-f(x)|\,d\mu(z).\label{boom}
\end{equation}

Then recall that \eqref{P2} admits a reformulation in terms of a pointwise estimate
(see \cite{H}, \cite{HaK}, Theorem 3.2), in particular it implies: 
\begin{equation}
|f(x)-f(z)|\le Cd(x,z)(h(x) + h(z)),\label{P3}
\end{equation} with
$$
h(x)=\left({\cal M}\left(|\nabla
f|^2\right)(x)\right)^{1/2}
$$
  for all $f$ with $f,\nabla f$ locally square integrable and $x,z\in M$.

Then, plugging (\ref{P3}) and \eqref{gradgub}  into (\ref{boom}),
$$|\nabla e^{-t\Delta}f(y)|
\le\frac{C}{V(x,\sqrt{t})}
\int_M\frac{d(x,z)}{\sqrt{t}}\exp\left(-\frac{d^2(y,z)}{Ct}\right)(h(x) + h(z))\,d\mu(z),
$$
and since $x,y\in B$, $\frac{d(x,z)}{\sqrt{t}}$ is comparable up to an additive constant with
$\frac{d(y,z)}{\sqrt{t}}$, therefore 
\begin{eqnarray*}
|\nabla e^{-t\Delta}f(y)|
&\le&\frac{C}{V(x,\sqrt{t})}
\int_M\exp\left(-\frac{d^2(y,z)}{C't}\right)(h(x) + h(z))\,d\mu(z)\\
&\le& C h(x)+\frac{C}{V(x,\sqrt{t})}
\int_M\exp\left(-\frac{d^2(y,z)}{C't}\right)h(z)
\,d\mu(z)\\
&\le& Ch(x)+C({\cal M}h)(x),
\end{eqnarray*}
again by
\cite{DR}, Proposition 2.4.
Now, a result of Coifman and Rochberg \cite{CR}, which is extendable to
spaces of homogeneous type, states that for any $g$ for which ${\cal M} g
< \infty$  a.e. and for any positive $\gamma < 1,$ the weight $w
\equiv ({\cal M} g)^{\gamma}$ belongs to the Muckenhoupt class $A_1$, with
$A_1$ constant depending on $\gamma$ (but not on $g$).  Thus,
${\cal M} w \leq C_{\gamma} w$ a.e., so that, in particular,
$$ {\cal M}\left[\left({\cal M}g\right)^{1/2}\right] \leq 
C_{1/2}\left({\cal M}g\right)^{1/2},$$
almost everywhere.

Applying this with $g=|\nabla f|^2$ in the right-hand side of the above inequality,
one obtains \eqref{domna}.

For the sake of completeness, we are now going to give a lemma which clarifies the relationship between
the pointwise  gradient bound \eqref{gradgub}, an integral version of it, the domination condition \eqref{conj}, and the
maximal estimate
\eqref{domna}. This lemma also offers a less direct, but more elementary approach to the implication from \eqref{gradgub}
to \eqref{domna}.

\begin{lem}\label{Steve}
In presence of $(FK)$, the following three properties are equivalent:

(i)  The  pointwise heat kernel gradient bound \eqref{gradgub}.

(ii) The weighted $L^2$ heat kernel gradient bound
\begin{equation}\label{intwry}
\int_{M} |\nabla_x
p_t(x,y)|^2e^{\alpha\frac{d^2(x,y)}{t}}\,d\mu(y)\le   \frac{C}{tV(x,\sqrt{t})},
\end{equation}
for $\alpha>0$ small enough, all $t>0$, $x\in M$.

(iii) The domination condition \eqref{conj}
\begin{equation*}
|\nabla e^{-t\Delta} f|^2\le C e^{-C't\Delta}(|\nabla f|^2), 
\end{equation*}
for some $C,C'>0$ and $f$ with $f, \nabla f$ square integrable and all $t>0$.

Any of these conditions implies the maximal estimate \eqref{domna}.

\end{lem}
 
\bigskip

\noindent{\bf Proof:} 
It is easy to see by integrating \eqref{gradgub} and using the doubling property that $(i)$ implies $(ii)$.
The converse  follows from an argument somewhat similar to  \cite{CDfull},  p.14.
Write
$$\nabla_x\,
p_{2t}(x,y) =  \int_M  \nabla_x\,
p_t(x,z) p_t(z,y) \, d\mu(z).$$
Using $(ii)$ and Cauchy-Schwarz inequality, 
$$|\nabla_x\,
p_{2t}(x,y)|^2 \le  \frac{C}{t\, V(x,\sqrt{t})} \int_M |p_t(z,y)|^2 
e^{-\alpha\frac{d^2(x,z)}{t}}\,d\mu(z)
$$
for all $\alpha>0$ small enough. 
According to $(DU\!E)$,
$$
|p_t(z,y)|^2  \le \frac{C}{V^2(y,\sqrt t) } e^{-c\frac{d^2(z,y)}{t}}
$$
for some $c>0$. For $\beta$ small enough we have 
 $$
e^{-\alpha\frac{d^2(x,z)}{t}}    e^{-c\frac{d^2(z,y)}{t}} \le  
e^{-\beta\frac{d^2(x,y)}{t}}e^{-c\frac{d^2(z,y)}{2t}},
$$
thus
$$|\nabla_x\,
p_{2t}(x,y)|^2 \le  \frac{C}{t\, V(x,\sqrt{t})V^2(y,\sqrt t)} e^{-\beta\frac{d^2(x,y)}{t}}\int_M
e^{-c\frac{d^2(z,y)}{2t}}   \,d\mu(z).
$$
 Using $(D)$, it is easy to check that the quantity
$$\frac{1}{V(y,\sqrt t)} \int_M
e^{-c\frac{d^2(z,y)}{2t}}   \,d\mu(z)$$
 is uniformly bounded,
and $(i)$ follows readily.

 Assume next that $(i)$ holds and let us   prove  $(iii)$.
Recall that we may use freely $(LY)$ and \eqref{P2}.  

Fix  $x\in M$ and $B \equiv B(x,\sqrt{t}).$ Set $C_1=4B$ and $C_j=2^{j+1}B \setminus 2^jB$ for
$j\ge 2$. Using again (\ref{stoco}), write 
\begin{eqnarray*}|\nabla e^{-t\Delta}f(x)|&=&|\nabla 
e^{-t\Delta}(f-f_{4B})(x)|\\
&=&\left|\int_M\nabla_x\, p_t(x,y)(f(y)-f_{4B})\,d\mu(y)\right|\\
&\le&\sum_{j\ge
1}\int_{C_j}|\nabla_x\, p_t(x,y)| |f(y)-f_{4B}|\,d\mu(y).
\end{eqnarray*}
It follows from the lower bound in $(LY)$ that, for  all $\varepsilon>0$, there exist $c_\epsilon,
C_\varepsilon>0$ independent of $x \in M$ and $t>0$ such  that
\begin{equation}\frac{c_\varepsilon}{V(x,\sqrt{t})}\exp(- \varepsilon 4^j) \leq
p_{C_\varepsilon t}(x,y)
\label{lbj}\end{equation}
for all $ y \in 2^j B$, $t>0$.

Let us treat the first term in the above sum, that is when $j=1$. According to  \eqref{gradgub},
$$\int_{4B}|\nabla_x\, p_t(x,y)| |f(y)-f_{4B}|\,d\mu(y)\le 
\frac{C}{\sqrt{t}V(x,\sqrt{t})}\int_{4B}
|f(y)-f_{4B}|\,d\mu(y).$$
By (\ref{P2}),
$$\int_{4B}
|f(y)-f_{4B}|\,d\mu(y)\le C\sqrt{t}\sqrt{V(x,\sqrt{t})}\left(\int_{4B}
|\nabla f(y)|^2\,d\mu(y)\right)^{1/2},$$
hence
$$\int_{4B}|\nabla_x\, p_t(x,y)| |f(y)-f_{4B}|\,d\mu(y)\le 
C\left(\frac{1}{V(x,\sqrt{t})}\int_{4B}
|\nabla f(y)|^2\,d\mu(y)\right)^{1/2},$$
and by (\ref{lbj}),
$$\int_{4B}|\nabla_x\, p_t(x,y)| |f(y)-f_{4B}|\,d\mu(y)\le C'_\varepsilon \left( \int_M
p_{C_\varepsilon t}(x,y)|\nabla f(y)|^2
\,d\mu(y) \right)^{1/2}.
$$
Now for the other terms in the sum.
By \eqref{gradgub} again,
\begin{eqnarray*}
\int_{C_j}|\nabla_x\, p_t(x,y)| |f(y)-f_{4B}|\,d\mu(y)&\le&
\frac{C}{\sqrt{t}V(x,\sqrt{t})}
\int_{C_j}\exp\left(-\frac{d^2 (x,y)}{Ct}\right)
|f(y)-f_{4B}|\,d\mu(y)\\
&\le&
\frac{C\exp(-c4^j)}{\sqrt{t}V(x,\sqrt{t})}
\int_{C_j}
|f(y)-f_{4B}|\,d\mu(y).
\end{eqnarray*}
Then
\begin{eqnarray*}
\int_{C_j}
|f(y)-f_{4B}|\,d\mu(y)
&\le&
\int_{2^{j+1}B}\left(|f(y)-f_{2^jB}|
+\sum_{\ell=2}^{j}|f_{2^\ell B}-f_{2^{\ell+1}B}|\right)\,d\mu(y)\\
&=&
\int_{2^{j+1}B}|f(y)-f_{2^jB}|\,d\mu(y)+
\sum_{\ell=2}^{j}V(x,2^j\sqrt{t})|f_{2^\ell B}-f_{2^{\ell +1}B}|.
\end{eqnarray*}
Again by (\ref{P2}),
$$\int_{2^{j+1}B}|f(y)-f_{2^jB}|\,d\mu(y)\le 
C\,2^j \sqrt{t}\sqrt{V(x,2^j\sqrt{t})}\left(
\int_{2^{j+1}B}|\nabla
f(y)|^2 \,d\mu(y) \right)^{1/2},$$
and by (\ref{lbj})
$$\int_{2^{j+1}B}|\nabla
f(y)|^2 \,d\mu(y)\le C_\varepsilon\,V(x,\sqrt{t})\exp( \varepsilon 4^j)
\int_Mp_{C_\varepsilon t}(x,y)|\nabla
f(y)|^2\,d\mu(y),$$
thus 
$$ \int_{2^{j+1}B}|f(y)-f_{2^jB}|\,d\mu(y)\le  C_\varepsilon\,2^j \sqrt{t}V(x,2^j\sqrt{t})\,
\exp(\varepsilon4^j)\left(\int_M p_{C_\varepsilon t}(x,y) |\nabla
f(y)|^2 \,d\mu(y) \right)^{1/2}.$$
Similarly,
\begin{eqnarray*}
|f_{2^\ell B}-f_{2^{\ell +1}B}|&\le&
\frac{1}{V(x,2^\ell\sqrt{t})}\int_{2^{\ell+1}B}|f(y)-f_{2^{\ell 
+1}B}|\,d\mu(y)\\
&\le&
C\,2^\ell\sqrt{t}\left(\frac{1}{V(x,2^{\ell}\sqrt{t})}
\int_{2^{\ell+1}B}|\nabla
f(y)|^2 \,d\mu(y) \right)^{1/2},\end{eqnarray*}
thus by (\ref{lbj}) again,
$$|f_{2^\ell B}-f_{2^{\ell +1}B}|\leq C_\varepsilon \, 2^\ell\sqrt{t} \, 
\exp(\varepsilon 4^{\ell+1})\left(\int_M p_{C_\varepsilon t}(x,y) |\nabla
f(y)|^2 \,d\mu(y) \right)^{1/2}.
$$
Hence
\begin{eqnarray*}
&&\int_{C_j}
|f(y)-f_{4B}|\,d\mu(y)\\
&\le&
C_\varepsilon\,\sqrt{t}V(x,2^j\sqrt{t})\left(2^j+
\sum_{\ell=2}^{j}2^\ell\right)\exp(\varepsilon 4^j)\left(\int_M p_{C_\varepsilon t}(x,y) |\nabla
f(y)|^2 \,d\mu(y) \right)^{1/2}\\
&\le&
C'_\varepsilon\,2^j\sqrt{t}V(x,2^j\sqrt{t})\exp(\varepsilon 4^j)\left(\int_M p_{C_\varepsilon
t}(x,y) |\nabla f(y)|^2 \,d\mu(y) \right)^{1/2}.
\end{eqnarray*}
Gathering the above estimates and using the doubling property, one obtains
\begin{eqnarray*}&&\int_{C_j}|\nabla_x\, p_t(x,y)| 
|f(y)-f_{4B}|\,d\mu(y)\\
&\le& 
C\,2^j\exp(-(c-\varepsilon)4^j)\frac{V(x,2^j\sqrt{t})}{V(x,\sqrt{t})}\left(\int_M p_{C_\varepsilon
t}(x,y) |\nabla f(y)|^2 \,d\mu(y) \right)^{1/2}\\
&\le& C'\,2^j\exp(-(c-\varepsilon)4^j)2^{\nu j}\left(\int_M p_{C_\varepsilon t}(x,y) |\nabla
f(y)|^2 \,d\mu(y) \right)^{1/2}.
\end{eqnarray*}
Choosing $\varepsilon<c$, 
since then $\sum_{j\ge 2}\,2^j\exp(-(c-\varepsilon)4^j)2^{\nu j}<\infty$,
$$|\nabla e^{-t\Delta}f(x)|\le C \, \left(\int_M p_{C_\varepsilon t}(x,y) |\nabla
f(y)|^2 \,d\mu(y) \right)^{1/2},$$
and $(iii)$ is proved.

The converse, that is the implication from $(iii)$ to $(i)$, again follows from a variant of the argument in \cite{CDfull},
p.14. Using again 
$$\nabla_x\,
p_{2t}(x,y) =  \int_M  \nabla_x\,
p_t(x,z) p_t(z,y) \, d\mu(z)$$ and 
$(iii)$, we have 
$$|\nabla_x\,
p_{2t}(x,y)|^2\le  C \int_M  
p_{C't}(x,z) |\nabla_z\, p_t(z,y)|^2 \, d\mu(z). $$
Invoke a weighted $L^2$ estimate for $\nabla_x\,p_t(.,y)$ proved in \cite{G1} (see also \cite{CD}, Lemma 2.4) and
valid under 
$(D)$ and $(DU\!E)$: for some $\gamma>0$ and all $y\in M, t>0$, 
\begin{equation}\label{gradx}
\int_M | \nabla_x\, p_t(x,y)|^2  e^{\gamma\frac{ d^2(x,y)}{t}} \, d\mu(x) \le \frac C{ t
\,V(y,\sqrt{t})}
\end{equation}
(note that contrary to $(ii)$, the integration here is with respect to $x$, see the Remark below).
This yields
$$|\nabla_x\,
p_{2t}(x,y)|^2\le  \frac{C }{t\, V(y,\sqrt t)} \left(\sup_{z\in M} 
p_{C't}(x,z) e^{-\gamma\frac{d^2(z,y)}{t}}\right) .$$
Using $(U\!E)$, the above supremum can easily be controlled by
$$  \frac{C e^{-\beta\frac{d^2(x,y)}{t}}}{ V(x,\sqrt t)}
$$
for   $\beta>0$ small enough, and $(i)$ follows.

\bigskip

It remains to deduce the maximal estimate \eqref{domna} from one of the other equivalent conditions. Assume $(iii)$, that is 
$$
|\nabla e^{-t\Delta} f(y)|^2\le C e^{-C't\Delta}(|\nabla f|^2)(y), 
$$
then, since  by $(LY)$,
$$e^{-ct\Delta}(|\nabla f|^2)(y)\le C e^{-C't\Delta}(|\nabla f|^2)(x)$$
as soon as $d(x,y)\le \sqrt{t}$,
one obtains 
$$
|\nabla e^{-t\Delta} f(y)|^2\le C e^{-C't\Delta}(|\nabla f|^2)(x). $$
On the other hand,
$$e^{-C't\Delta}(|\nabla f|^2)(x)\le C {\cal M}(|\nabla f|^2)(x)$$
by  $(U\!E)$ and \cite{DR}, Proposition 2.4.  This readily yields \eqref{domna}.

\bigskip

\noindent{\bf Remarks:}

-Although  not necessary for the main argument developed in this section, we have recorded the equivalence
between
$(i)$ and
$(ii)$ to point out a difference with the  case 
$1<p<2$ (see \cite{CD}).  In that case, the crucial ingredient in the proof was the weighted estimate  of the gradient
\eqref{gradx}, where integration and differentiation are taken with respect to the same variable. As  we already said,
\eqref{gradx} follows  from the pointwise estimate of the heat kernel only. 
   In contrast, the estimate
 $(ii)$, with integration and differentiation with respect to different variables, requires in addition the pointwise estimate
$(i)$ of the gradient. All this also explains at a technical level why more assumptions are needed in Theorem \ref{maincor}
than in Theorem \ref{CD}, and also why
\eqref{gradgub} holds 
 on  manifolds satisfying $(FK)$ and $$|\nabla_xp_t(x,y)|\le C |\nabla_y
p_t(x,y)|$$  for all $t>0,\,x,y\in M$  (\cite{G1}, Theorem 1.3).

-A sufficient condition for $(iii)$   in terms of Ricci curvature follows from \cite{LXD}, Theorem 3.2, (3.7).
According to the
above,  $(R_p)$ holds for all $p\in (1,\infty)$ on manifolds satisfying this condition plus $(FK)$.

- The  $L^2$ version of Poincar\'e inequalities
$(P)$, which follows from  the assumptions of the Lemma and $(i)$, is used in the above proof. If instead one has the stronger
$L^1$ version
$$
\int_{B}|f-f_{B}|\,d\mu\le Cr\int_{B}
|\nabla f|\,d\mu,
$$
for any ball $B$ in $(M,d)$ with radius $r$  and any $f$ with $f, \nabla f$ locally integrable on
$B$ (which is the case for instance on Lie groups of polynomial volume growth, manifolds with
non-negative Ricci curvature and co-compact coverings with polynomial growth deck transformation
group),  together $(D)$, $(DU\!E)$, then one  can show in the same way the equivalence of \eqref{gradgub} with the estimates
$$
\sup_{y\in B}|\nabla e^{-r^2(B)\Delta}f(y)|\le C'\inf_{x\in B}{\cal M} (|\nabla f|)(x),
$$
and
$$
|\nabla e^{-t\Delta} f|\le C e^{-t\Delta}(|\nabla f|), 
$$
which are stronger than their $L^2$ counterparts.

\subsection{Proofs of some other results}\setcounter{equation}{0}\label{other}

In this section, $M$  satisfies $(D)$ and $(DU\!E)$.

\bigskip

\noindent{\bf Proof of Proposition \ref{proother}:} First assume $(G_p)$. 
Write
$$\||\nabla_x\,p_{2t}(.,y)|\|_p=\| |\nabla e^{-t\Delta}(p_{t}(.,y))|\|_p \le \frac{C_p}{ \sqrt{t}}\|
p_t(.,y)\|_p.
$$
The estimate $(U\!E)$ easily yields
$$\norm{p_t(.,y)}_p\le \frac{C}{ \left[ V(y,\sqrt{t})\right]^{1-\frac{1}{ p}}}.$$
This implies
$$\norm{|\nabla_x\,p_{2t}(.,y)|}_p\le  \frac{C_p}{ \sqrt{t}\left[ V(y,\sqrt{t})\right]^{1-\frac{1}{
p}}}.$$

Conversely, assume for a $p \in (2,p_0)$ for all $y\in M$ and $t>0$.
 $$\norm{|\nabla_x\, p_t(.,y)|}_p\le  \frac{C_p}{ \sqrt{t}\left[ V(y,\sqrt{t})\right]^{1-\frac{1}{
p}}}.$$
 We shall prove $(G_q)$ for any $2<q<p$. 
Using \eqref{gradx} and
interpolating with the above unweighted $L^p$ estimate, we have
$$
\int_M | \nabla_x\, p_t(x,y)|^q  e^{\gamma'\frac{ d^2(x,y)}{t}} \, d\mu(x) \le \frac C{ t^{q/2}
\left[V(y,\sqrt{t})\right]^{q-1}}.
$$
Now let $f \in L^q(M,\mu)$. Estimate 
$|\nabla e^{-t\Delta}f(x)|$ by
$$
\int_M |\nabla_x\, p_t(x,y)|  {e^{\frac{\gamma' d^2(x,y)}{qt}}}{\left[V(y,\sqrt{t})\right]^{1/q'}}
|f(y)|
\frac {e^{-\frac{ \gamma'd^2(x,y)}{qt}}}{\left[V(y,\sqrt{t})\right]^{1/q'}}\, d\mu(y)
$$
which, by H\"older inequality, is controlled by
$$
\left(\int_M | \nabla_x\,p_t(x,y)|^q  {e^{\gamma'\frac{
d^2(x,y)}{t}}}{\left[V(y,\sqrt{t})\right]^{q/q'}} |f(y)|^q\, d\mu(y) \right)^{1/q}
\left(\int_M \frac {e^{-\frac{q'\gamma' d^2(x,y)}{qt}}}{V(y,\sqrt{t})}\,
d\mu(y)\right)^{1/q'}.
$$ 
Now, it follows easily from the doubling property that   the second integral is bounded uniformly in $t,x$. 
Integrating with respect to $x$, by Fubini's theorem and the weighted $L^q$-estimate, we obtain
$$
\int_M |\nabla e^{-t\Delta}f(x)|^q \, d\mu(x) \le \frac {C}{t^{q/2}} \int_M |f(y)|^q\, d\mu(y)
$$
as desired. 

\bigskip

\noindent{\bf Proof of Theorem \ref{maincor}:} We have already observed that the hypotheses in the
statement imply $(D)$ and $(P)$. It remains to obtain  $(G_p)$
for all $p\in (1,\infty)$. In view of the previous result, it is enough to prove \eqref{gradlp}.
But this follows by interpolating  $(G)$ with the  
$L^2$ bound
$$\norm{|\nabla_x\,p_t(.,y)|}_2\le  \frac{C_2}{ \sqrt{t}\left[ V(y,\sqrt{t})\right]^{\frac{1}{ 2}}}.$$
The latter follows from $(DU\!E)$:
\begin{eqnarray*}
\norm{|\nabla_x\,p_t(.,y)|}_2^2&=&(\Delta p_t(.,y),p_t(.,y))\\
&\le&\norm{\Delta p_t(.,y)}_2\norm{p_t(.,y)}_2.
\end{eqnarray*} 
Now
$$\norm{\Delta p_t(.,y)}_2=\norm{\Delta e^{-(t/2)\Delta}p_{t/2}(.,y)}_2\le \frac{C}{t}\norm{p_{t/2}(.,y)}_2,$$
by analyticity of the heat semigroup on $L^2$, and $\norm{p_{t/2}(.,y)}_2^2=p_t(y,y)\le\frac{C}{V(y,\sqrt{t})}$ by $(DU\!E)$.
The claim is proved.

\section{Localization}\setcounter{equation}{0}\label{loc}

This section is devoted to the proof of  Theorem \ref{mainp>2loc}. Theorem \ref{maincorloc}
can then be deduced as in the global case and we leave details to the reader.

 For $a>0$, we have
$$\nabla(\Delta+a)^{-1/2}=c\int_0^{\infty}e^{-at}\nabla e^{-t\Delta}{
\frac{dt}{ \sqrt{t}}}.$$

Let $v$ be a $C^\infty$  function on $[0,\infty)$ with $0\le v\le 1$,  which equals $1$ 
on
$[0, 3/4]$ and  vanishes on
$t\ge 1$.

If $a>\alpha$,  Minkowski's integral  inequality implies
that 
$$ \left\|\int_{0}^{\infty}(1-v(t))e^{-at}\nabla e^{-t\Delta}f{
\frac{dt}{ \sqrt{t}}}\right\|_{p} \le C \int_{3/4}^\infty e^{(\alpha-a)t} {
\frac{dt}{ \sqrt{t}}} \|f\|_{p}\le C'\|f\|_{p}.
$$
It is enough to prove the $L^p$ boundedness of the sublinear operator ${\tilde T}$ defined by  
$${\tilde T}f=\left|\int_0^{\infty} v(t)\nabla e^{-t(\Delta +a)}f
\frac{dt}{ \sqrt{t}}\right|,$$
 Without loss of generality, we may and do replace in what follows $\Delta +a$ by $\Delta$ as the value of 
$a$ 
 plays no further role.

For this purpose, we begin with the localization technique of
\cite{DR} as in
\cite{CD}.

Before we start,  observe that, as a
consequence of
$(E)$ and
$(P_{loc})$, $p_t(x,y)$ satisfies the estimates $(LY)$ for small times.

Let $(x_j)_{j\in J}$ be a maximal $1$-separated subset of $M$: the collection of
 balls $B^j=B(x_j,1)$,
 $j\in J$, covers $M$, whereas  the balls $B(x_j, 1/2)$ are pairwise disjoint. 
It follows from
the local doubling property that there exists $N\in \nat^*$ such that 
every $x\in M$ is contained in at
most $N$ balls $4B^j=B(x_j,4)$.

Consider a $C^\infty$ partition of unity $\varphi_j$, $j\in J$, such that 
$\varphi_j \ge 0$ and is supported in  $B^j$. Let  $\chi_j$ be the characteristic function of the
ball $4B^j$. For $f\in \ccc_0^\infty(M)$ and $x\in M$, write
$${\tilde T}f(x)\le \sum_{j}\chi_j{\tilde
T}(f\varphi_j)(x)+\sum_{j}(1-\chi_j){\tilde T}(f\varphi_j)(x)= I+II.  $$

Let us first treat the term $II$. 
The first observation is that, from the finite overlap property of the balls $B^j$, we have 
$$
 \sum_j | (1-\chi_j)(x) \varphi_j(y)| \le N  \chi_{d(x,y)\ge 3}.
$$
Hence
\begin{align*}
II& \le  \int_0^1\int_M  |\nabla_xp_t(x,y)|  \sum_j | (1-\chi_j)(x) \varphi_j(y)| |f(y)|\, d\mu(y) \frac
{dt}{\sqrt t}\\
&\le N \int_0^1\int_{d(x,y)\ge 3}  |\nabla_xp_t(x,y)|     |f(y)|\, d\mu(y) \frac
{dt}{\sqrt t}.
\end{align*}
From there, we follow the argument in Section \ref{other} by inserting the Gaussian terms and using
H\"older inequality to estimate the integral on ${d(x,y)\ge 3}$.  Since $(LY)$ for small times applies,
we have
$$
\int_{d(x,y)\ge 3} \frac {e^{-\frac{p'\gamma d^2(x,y)}{pt}}}{V(y,\sqrt{t})}\,
d\mu(y) \le C e^{-\frac c t } 
$$
for some constants $C,c>0$ independent of $t$ and $x$.
Therefore, 
$$II \le C \int_0^1\left(\int_M |\sqrt t\,  \nabla_xp_t(x,y) |^p {e^{\gamma\frac{
d^2(x,y)}{t}}}{\left[V(y,\sqrt{t})\right]^{p/p'}} |f(y)|^p\, d\mu(y) \right)^{1/p} \, 
\frac {e^{-\frac
c t } }{t}\, dt.
$$
It follows from $(G_{p_0}^{loc})$ and the argument in Section \ref{other}  for small times (which is valid under the
exponential growth assumption and the Gaussian upper bound for $p_t(x,y)$ for small times, see \cite{CD}) that for some
$\gamma>0$ and all $t\in (0,1)$ and $y\in M$,  
$$
\int_M | \nabla_xp_t(x,y)|^p  e^{\gamma'\frac{ d^2(x,y)}{t}} \, d\mu(x) \le \frac C{ t^{p/2}
\left[V(y,\sqrt{t})\right]^{p-1}}.
$$
Since the measure $\frac {e^{-\frac c
t } }{t}\, dt$ has finite mass, one can use Jensen's inequality with respect to $t$, Fubini's theorem
and the weighted
$L^p$ estimate above to conclude that $\int_M | II |^p \, d\mu(x) $ is bounded by $C\|f\|_p^p$.

We now turn to $I\,$ which is the main term. The uniform
overlap of the balls
$4B^j$ implies 
$$
\sum_j \|g\chi_j\|_r^r \le C \|g\|_r^r
$$for all $g$ and $1\le r\le \infty$. Hence 
$$
\int_M |g(x)|\bigg|\sum_{j}\chi_j{\tilde
T}(f\varphi_j)(x)\bigg|\, d\mu(x) \le 
C \sum_{j} \|f\varphi_j\|_p\|g\chi_j \|_{p'} 
\le C
\|f\|_p\|g\|_q,
$$
provided we show that 
$$
\|\chi_j |{\tilde T}(f\varphi_j)|\|_p \le C\|f\varphi_j\|_p
$$
with a bound uniform in $j$.  In other words, we want to show that $\tilde T$ maps $L^p(B^j)$
into 
$L^p(4B^j)$. To this end we apply Theorem \ref{lploc} with $E_1=B^j$ and $E_2=4B^j$ since $4B^j$
has the doubling property by  the following lemma, which is implicit in 
\cite{LHQ}, and whose proof we postpone until the end of this section.

\begin{lem}\label{doublingballs} The balls $4B^j$ equipped with the induced distance and measure satisfy
the doubling property
$(D)$ and the doubling constant may be chosen independently of $j$. More precisely, there is a constant
$C\ge 0$ such that for all $j\in J$, 
\begin{equation} \label{Dj}
\mu(B(x,2r) \cap 4B^j)\le C\, \mu(B(x,r)\cap 4B^j),\ \forall\,x\in 4B^j,\,r>0,
\end{equation}
and also
\begin{equation} \label{Djcap}
\mu(B(x,r))\le C\,\mu(B(x,r)\cap 4B^j),\ \forall\,x\in 4B^j,\,0<r\le 8.
\end{equation}
\end{lem}

  Define the local maximal
function on $M$ by
$${\cal M}^{loc}f(x)=\sup_{B\ni x,  r(B) \le 32}\frac{1}{\mu(B)}\int_B|f|\,d\mu,$$ for $x\in M$ and $f$ locally integrable on
$M$.
 By local doubling, ${\cal M}^{loc}$  is bounded on all $L^p(M)$, $1<p\le \infty$. 
We have the following estimates for balls $B$ of $M$ centered in $4B^j$ and with radii less than 8, $x
\in B\cap 4B^j$ and functions $f$ in $L^2(M)$ with support in $B^j$ :

1) There is an integer $n$ depending only on the condition $(E)$ such that  
\begin{equation}\label{pointwiserieszloc}
\frac 1{\mu(B\cap 4B^j)}\int_{B\cap 4B^j} |{\tilde T}(I-e^{-r^2\Delta})^nf|^2\, d\mu \le C    {\cal
M}^{loc}(|f|^2)(x).
\end{equation} 

2)  If $2<p<p_0$, then for $I-A_r= (I-e^{-r^2\Delta})^n$ 
\begin{equation}\label{maxloc}
\left( \frac 1{\mu(B\cap 4B^j)}\int_{B\cap 4B^j}   |{\tilde T} A_rf|^p\, d\mu
\right)^{1/p}  
\le  C
\big({\cal M}_{4B^j}
(|{\tilde T} f|^2)+ (Sf)^2\big)^{1/2}(x)
\end{equation}
with 
\begin{equation}
(Sf)^2= {\cal
M}^{loc}(|{\tilde T} f|^2 \chi_{(M\setminus 4B^j)}) +  {\cal
M}^{loc}(| h|)^2 + {\cal M}_{4B^j}
(|f|^2)
\end{equation}
and  $h = \int_0^{\infty} v(t)e^{-t\Delta }f{
\frac{dt}{ \sqrt{t}}}$.

Admitting these inequalities, it remains to see that $\|Sf\|_{L^p(4B^j)} \le C\|f\|_{
L^p(B^j)} $. By the study of $II$, we have that 
$$
\|\big({\cal
M}^{loc}(|{\tilde T} f|^2 \chi_{(M\setminus 4B^j)}\big)^{1/2}\|_{ L^p(4B^j)} \le 
C\||{\tilde T} f| \chi_{(M\setminus 4B^j)}\|_{ L^p(M)} \le   C\|f\|_{ L^p(B^j)}.
$$
Next, by definition of $h$ and the contraction property of the heat semigroup, $\|h\|_{L^p(M)} \le 
c \|f\|_{ L^p(B^j)}$, hence 
$$
\|{\cal
M}^{loc}(|h| )\|_{ L^p(4B^j)} \le 
   C\|f\|_{ L^p(B^j)}.
$$
Finally, we conclude the argument by invoking the boundedness of  ${\cal M}_{4B^j}$  
  on
$L^{p/2}(4B^j)$ to bound $\|\big({\cal M}_{4B^j}
(|f|^2)\big)^{1/2}\|_p$.
\bigskip

\noindent{\bf  Proof of (\ref{pointwiserieszloc}):}  Take a ball $B$ centered in $4B^j$ with $r=r(B)<8$
and  let $x$ be  any point in $B\cap 4B^j$. By Lemma \ref{doublingballs} we have 
$\mu(B) \le C\mu(B\cap 4B^j)$ for some $C$ independent of $B$ and $j$.  Hence, 
$$\frac 1{\mu(B\cap 4B^j)}\int_{B\cap 4B^j} |{\tilde T}
(I-e^{-r^2\Delta})^nf|^2\, d\mu \le \frac C{\mu(B)} \int_{B} |{\tilde T}
(I-e^{-r^2\Delta})^nf|^2\, d\mu 
$$ and we follow the calculations of Section \ref{Riesz} with $\tilde T$ replacing $\nabla
\Delta^{-1/2}$.  Introduce $i_r$ the integer defined by $2^{i_r} r < 8 \le 2^{i_r+1} r$.  Denote by
$C_i$ the ring $2^{i+1}B\setminus 2^iB$ if $i\ge 2$ and $C_1=4B$. Decompose $f$ as $f_1+f_2+f_3+\ldots
+ f_{i_r}$ with 
$f_i=f\chi_{C_i}$. The decomposition stops since $f$ is supported in $B^j$ and 
 $4B^j \subset 2^i B$ when $i>i_r$. 
By Minkowski's
inequality we have that 
$$
\| |\tilde T(I-e^{-r^2\Delta})^nf|\|_{L^2(B)} \le \sum_{i= 1}^{i_r} 
\| |\tilde T
(I-e^{-r^2\Delta})^nf_i|\|_{L^2(B)}.
$$
For $i=1$ we use the $L^2$ boundedness of $\tilde T
(I-e^{-r^2\Delta})^n$:
$$
\| |\tilde T
(I-e^{-r^2\Delta})^nf_1|\|_{L^2(B)} \le \| f \|_{L^2(4B)} \le \mu(4B)^{1/2}
 \big({\cal M}^{loc}(|f|^2)(x)\big)^{1/2},
$$
For $i\ge 2$ we use the definition of $\tilde T$:
\begin{eqnarray*}\tilde T
(I-e^{-r^2\Delta})^nf&=&\left|\int_0^\infty v(t) \nabla e^{-t\Delta}f
(I-e^{-r^2\Delta})^n\,\frac{dt}{\sqrt t}\right|\\ &=&\left|\int_0^\infty  \tilde g_{r}(t) \nabla e^{-t\Delta}f
\, {dt}\right|\\
\end{eqnarray*}
where using the usual notation for the binomial coefficient, 
$$\tilde g_{r}(t)=\sum_{k=0}^n \binom n k (-1)^k \frac{ v(t-kr^2)}
{\sqrt{t-kr^2}}\chi_{\{t>kr^2\}}.$$  
By Minkowski's integral inequality and the Gaffney estimates (\ref{Gaffney}) using the support of $f_i$, we
have that 
$$
\| |\tilde T
(I-e^{-r^2\Delta})^n f_i|\|_{L^2(B)} \le C \int_0^\infty |\tilde g_{r}(t)| e^{-\frac
{\alpha' 4^ir^2} t}  \frac {dt} {\sqrt t} \ \| f\|_{L^2(C_i)}.
$$ 
The latter integral can be estimated as follows.  
 Elementary analysis yields the following estimates for $\tilde g_{r}(t)$:
$$
|\tilde g_{r}(t)| \le \frac C{\sqrt {t-kr^2}} \quad \textrm{if} \quad {kr^2<t \le(k+1)r^2 \le
(n+1)r^2},
$$
$$
|\tilde g_{r}(t)| \le Cr^{2n} \quad \textrm{if} \quad{ (n+1)r^2 < t},
$$
and $\tilde g_{r}(t)=0$ for $t \ge 1 + nr^2$. 
Hence 
$$
\int_0^\infty |\tilde g_{r}(t)| e^{-\frac
{\alpha' 4^ir^2} t}  \frac {dt} {\sqrt t} \le C \min \{4^{-in}, r^{2n}\}= C 8^{-2n} r^{2n}.
$$
Now, an easy consequence of local doubling and $r(2^i B) \le 8$ when $1\le i \le i_r$,  is that  
 $$\mu(2^{i+1} B) \le m(2^{i+1}) \mu(B)$$  with $m(\theta)=C(1+\theta)^\nu$,  $C$ and $\nu$ independent of
$B$ and
$j$. Therefore, as $C_i \subset 2^{i+1}B$, 
$$\| f\|_{L^2(C_i)}  \le \mu(2^{i+1}B)^{1/2}
 ({\cal M}^{loc}(|f|^2)(x))^{1/2} \le  \sqrt{m(2^{i+1})} \mu(B)^{1/2}  \big({\cal
M}^{loc}(|f|^2)(x)\big)^{1/2}.
$$
Choosing $2n>\nu/2$ and using the definition of $i_r$ and $r\le 8$, we obtain 
\begin{align*}
\| |\tilde T
(I-e^{-r^2\Delta})^n f|\|_{L^2(B)} &\le C \sum_{i= 1}^{i_r}  2^{i\nu/2}r^{2n} \mu(B)^{1/2} 
\big({\cal M}^{loc}(|f|^2)(x)\big)^{1/2}
\\
& \le C  \mu(B)^{1/2} 
\big({\cal M}^{loc}(|f|^2)(x)\big)^{1/2}.
\end{align*}

\bigskip

\noindent{\bf  Proof of (\ref{maxloc}):} We first establish  the analog of Lemma \ref{lemoffdiagl2lp}:
for every
$p\in (2,p_0)$, for every ball
$B$ centered in $4B^j$ with radius
$r< 8$ and every
$L^2$ function $g$ supported in $C_i=2^{i+1}B\setminus 2^iB$ if $i\ge 2$ or  $C_1=4B$ and every $k
\in
\{1,
\ldots, n\}$, where $n$ is chosen as above, one has
\begin{equation} \label{eqoffdiagl2lploc}
\left( \frac 1{\mu(B \cap 4B^j)}\int_{B\cap 4B^j}   |\nabla e^{-kr^2\Delta}g|^p\, d\mu \right)^{1/p}  
\le 
\frac {Ce^{- \alpha 4^i}}{r}\left(\frac 1 {\mu(2^{i+1}B)} \int_{C_i} |g|^2\, d\mu\right)^{1/2}.
\end{equation}
for some constants $C$ and $\alpha$ depending only on $(E)$, $(P_{loc})$,  $p$ and $p_0$.

The proof given in Section \ref{Riesz} can be copied \textit{in extenso} provided we make three remarks.
First, as in the previous argument since $\mu(B) \le C \mu(B\cap 4B^j)$ we may replace $B\cap 4B^j$ by
$B$ in the left-hand side of (\ref{eqoffdiagl2lploc}). Second,
as we already said, $p_t(x,y)$ satisfies the estimates $(LY)$ for, say,  $t=kr^2\le 64 n$.  Third, the polynomial
volume growth is replaced by an exponential volume growth but the estimates still carry out in this case thanks to the
Gaussian terms in the sums. Further details are left to the reader. 

Next, assume that $f\in L^2(B^j)$ and take $h = \int_0^{\infty} v(t)e^{-t\Delta }f{
\frac{dt}{ \sqrt{t}}}$. Since $\tilde T f= |\nabla h|$ and according to the first remark above, it is
enough to control
$$
 \left( \frac 1{\mu(B)}\int_B   |\nabla e^{-kr^2\Delta}h|^p\, d\mu \right)^{1/p}   .
$$
Write $\nabla e^{-kr^2\Delta}h= \nabla e^{-kr^2\Delta}(h-h_{4B})= \sum_{i\ge 1} \nabla
e^{-kr^2\Delta}g_i$ where $g_i=(h-h_{4B})\chi_{C_i}$.
Applying 
(\ref{eqoffdiagl2lploc}) to each $g_i$, we are reduced to estimating 
$\left(\frac 1 {\mu(2^{i+1}B)} \int_{C_i} |g_i|^2\, d\mu\right)^{1/2}$.
 We distinguish the two regimes $i\le
i_r$ and
$i >i_r$ where
$i_r$ is the largest integer satisfying  $2^i r <8$.  In the first regime, 
the argument in Section \ref{Riesz} using the local Poincar\'e inequalities for balls with radii not
exceeding 16 can be repeated and (\ref{eqci}) becomes
 \begin{equation}\label{eqciloc}
\left(\frac {1} {\mu(2^{i+1}B)} \int_{C_i}
|g_i|^2\,d\mu \right)^{1/2} \le  C \sum_{\ell=1}^i  (2^{\ell}r)
\left(\frac {1} {\mu(2^{\ell+1}B)} \int_{2^{\ell+1}B}
|\nabla h|^2\,d\mu \right)^{1/2}.
\end{equation} 
Write then
\begin{align*}
\frac {1} {\mu(2^{\ell+1}B)} \int_{2^{\ell+1}B}
|\nabla h|^2\,d\mu 
&\le \frac {1} {\mu(2^{\ell+1}B \cap 4B^j)} \int_{2^{\ell+1}B \cap 4B^j}
|\nabla h|^2\,d\mu \\
& \qquad\qquad\qquad+ \frac {1} {\mu(2^{\ell+1}B)} \int_{2^{\ell+1}B}
\chi_{(M\setminus 4B^j)}|\nabla h|^2\,d\mu
\\
&\le \big({\cal M}_{4B^j}
(|\nabla h|^2)(x)+  {\cal
M}^{loc}(|\nabla h|^2 \chi_{(M\setminus 4B^j)}\big)(x)
\end{align*}
where $x$ is any point in $B\cap 4B^j$.
Hence the contribution of the terms in the first regime does not exceed
$$ \sum_{1\le \ell \le i \le i_r} \frac {Ce^{- \alpha 4^i}}{r}
(2^{\ell}r) \big( ({\cal M}_{4B^j}
(|\nabla h|^2)(x)+  {\cal
M}^{loc}(|\nabla h|^2 \chi_{(M\setminus 4B^j)})(x)\big)^{1/2}.
$$
For the second regime, we proceed directly by
$$
\left(\frac 1 {\mu(2^{i+1}B)} \int_{C_i} |g_i|^2\, d\mu\right)^{1/2} 
\le
\left(\frac 1 {\mu(2^{i+1}B)} \int_{C_i} |h|^2\, d\mu\right)^{1/2} + |h_{4B}|.
$$
First $|h_{4B}|  \le {\cal
M}^{loc}(| h|)(x)$ since $B$ has radius less than 8. 
Next,  
$$
\int_{C_i} |h|^2\, d\mu \le \int_{M} |h|^2\, d\mu \le  C\int_M |f|^2\, d\mu= C\int_{4B^j} |f|^2\, d\mu
$$ 
since $f$ is supported in $B^j\subset 4B^j$. Also since $i>i_r$ we have $2^{i+1}B \supset 4B^j$ and 
$$
\left(\frac 1 {\mu(2^{i+1}B)} \int_{C_i} |h|^2\, d\mu\right)^{1/2} \le C \left(\frac 1 {\mu(4B^j)}
\int_{4B^j} |f|^2\, d\mu\right)^{1/2} \le C \big({\cal M}_{4B^j}
(|f|^2)(x)\big)^{1/2}.
$$
The contribution of the terms in the second regime is bounded above by
$$
\sum_{i>i_r} \frac {Ce^{- \alpha 4^i}}{r} \big({\cal
M}^{loc}(| h|)(x) + \big({\cal M}_{4B^j}
(|f|^2)(x)\big)^{1/2}\big)
$$
and it remains to recall that $1/r \le 2^i /8$ when $i>i_r$ to conclude the proof of 
(\ref{maxloc}). 

The  proof of Theorem \ref{mainp>2loc} is complete  provided we prove Lemma \ref{doublingballs} which we do
now.   

\bigskip

 We begin with the first inequality. Fix $j\in J$, $x\in 4B^j$, $r>0$. If $r\ge 8$, there is
nothing to prove. We assume $r<8$. There is a point $x_*$ such that  $B(x_*,
r/8) \subset 4B^j$ and $d(x,x_*)\le 3r/8$. Indeed, if $d(x,x_j)\le 3r/8$ then $B(x,
r/8) \subset 4B^j$, so that one can take $x_*=x$. Otherwise,  since $M$ is connected, there is a curve joining $x$
and $x_j$ whose length is smaller than, say, $d(x,x_j)+\frac{2r}{8}$. On this curve, one can choose $x_*$ such that $d(x,x_*)=
3r/8$, thus
$d(x_*,x_j)\le d(x,x_j)-\frac{r}{8}$.  The point $x_*$ satisfies the required properties.  Then one may write 
$$
\mu(B(x,2r) \cap 4B^j) \le V(x,2r) \le V(x_*, 19r/8) \le C V(x_*, r/8) \le C\mu(B(x,r)\cap 4B^j)
$$
where one uses the local doubling property for balls with radii not exceeding 19. 
This proves the first inequality of the lemma and the proof of the second one is contained in the
argument.

\section{$L^p$ Hodge decomposition for non-compact manifolds}\setcounter{equation}{0}\label{Hodge}

The Hodge decomposition, which associates to a form its exact, 
co-exact and harmonic parts, is
well-known to be bounded on $L^2$ on any complete Riemannian manifold 
(see for instance
\cite{DRh},  Theorem 24, p.165, and the recent survey \cite{Ca}).
The  question of the $L^p$ Hodge decomposition has been mainly examined on 
closed manifolds (see
\cite{Scott}) and or domains therein  (\cite{ISS}). On
  non-compact manifolds,  the
connection with the Riesz transform was established in the case of $1$-forms and an example was treated in \cite{St}.   The
unpublished manuscript \cite{LoDR} contains results for forms of all degrees in the case of 
 manifolds with positive bottom of the spectrum; see also \cite{LoDR}.  In the case of degree one, one
can deduce more general results from Theorems \ref{mainp>2} and 
\ref{pangap0}. We denote by
$L^pT^*M$  the usual $L^p$ space of $1$-forms.

\begin{theo}\label{h} Let $M$ be a complete non-compact Riemannian 
manifold satisfying either $(D)$, $(P)$ and $(G_{p_{0}})$ for some $p_0\in(2,\infty]$, or the assumptions of Theorem
${\rm \ref{pangap0}}$. Then the 
Hodge projector from  $1$-forms
onto   exact forms is bounded  on $L^pT^*M$, for all $p\in(q_0,p_0)$ where $q_0$ is the conjugate exponent to $p_0$.
\end{theo}

\noindent{\bf Proof:}
The projector on exact forms is $$d\Delta^{-1}\delta=d\Delta^{-1/2} 
(d\Delta^{-1/2})^*.$$
Now  $d\Delta^{-1/2}$ is bounded from $L^p(M,\mu)$ to $L^pT^*M$ for 
all $p\in(q_0,p_0)$ by the results in \cite{CD} and
Theorem \ref{mainp>2} or Theorem \ref{pangap0}, depending on the 
assumptions. By duality,  $(d\Delta^{-1/2})^*$ is bounded
  from $L^pT^*M$ to $L^p(M,\mu)$, for $p$ in the same range, and the claim follows.

\begin{cor}\label{hdec} Let $M$ be a complete non-compact Riemannian 
manifold satisfying either $(D)$, $(P)$ and $(G_{p_{0}})$ for some $p_0\in(2,\infty]$,   or the assumptions of Theorem
${\rm \ref{pangap0}}$. Let $p\in(q_0,p_0)$. Then the 
Hodge decomposition on $1$-forms extends to $L^pT^*M$ if one of the following additional assumptions holds:

a. $M$ is a surface.

b. There are no $L^p$ harmonic  $1$-forms.

c. The heat semigroup on $1$-forms,
$e^{-t\overrightarrow{\Delta}}$, is bounded on $L^pT^*M$ uniformly in $t>0$.

\end{cor}

\noindent{\bf Proof:} Let $p\in(q_0,p_0)$.  The projector on exact forms is bounded on 
$L^pT^*M$  by Theorem \ref{h}. In case
$a)$,  we  use an argument from
\cite{St}: in dimension $2$, the projector
   on co-exact forms is nothing but  $d\Delta^{-1}\delta$ conjugated 
by the Hodge star operator,
therefore
it is also bounded on $L^pT^*M$.
In  cases $b)$ and $c)$, the projection on harmonic forms extends boundedly
to
$L^pT^*M$:  In case $b)$, this projector is trivial by assumption, 
and in case $c)$, it is bounded
because it is the limit of $e^{-t\overrightarrow{\Delta}}$ as $t\to\infty$.
In all three
cases, two out of the three projectors in the  Hodge decomposition 
extend boundedly to $L^pT^*M$, therefore the
whole decomposition is an
$L^p$ decomposition.

\bigskip

Since,  on manifolds with
non-negative Ricci curvature, the Riesz transform is bounded on $L^p$ 
for all $p\in (1,\infty)$, and for all $t>0$ and $\omega\in
\ccc^\infty T^*M$,
$$|e^{-t\overrightarrow{\Delta}}\omega|\le  e^{-t\Delta}|\omega|,
$$    the proofs 
of Theorem
\ref{h} and Corollary \ref{hdec}
  show in particular that the Hodge decomposition extends to
$L^p$ on such manifolds. This fact seems known to experts, although 
we could not find it stated in the literature.

More generally, if the heat semigroup on forms  on $M$ is dominated 
by the heat semigroup on functions in the
  way discussed in Section \ref{assumption}:
\begin{equation}|e^{-t\overrightarrow{\Delta}}\omega|\le C 
e^{-ct\Delta}|\omega|,
\label{semidom}
\end{equation}
for some $C,c>0$ and all $t>0$ and $\omega\in
\ccc^\infty T^*M$, then both assumptions $b)$  (because there are no 
$L^p$ harmonic functions
  on a complete non-compact Riemannian manifold, see \cite{Ya}) and 
$c)$ are satisfied  for all $p\in(1,\infty)$.
It was proved in \cite{CDfull} (and it also follows from  Theorem 
\ref{maincor})  that
$(D)$,
$(DU\!E)$ and (\ref{semidom}) imply  the boundedness of the Riesz 
transform  for all $p\in(1,\infty)$.  Therefore, we can state
the following.

\begin{cor}\label{hdec1} Let $M$ be a complete non-compact Riemannian 
manifold satisfying
$(D)$, $(DU\!E)$ and ${\rm (\ref{semidom})}$. Then the Hodge 
decomposition on $1$-forms
extends to $L^pT^*M$ for all
$1<p<\infty$.\end{cor}

A similar statement could be formulated, with suitable assumptions, 
for manifolds where the bottom of the spectrum is
positive,
by using Theorem \ref{pangap}.

Finally, recall that the $L^p$ boundedness of the Hodge decomposition 
may be useful, together with other ingredients,
to show that, in certains situations, $e^{-t\overrightarrow{\Delta}}$ 
cannot be contractive on $L^pT^*M$ for all $p$
(see \cite{St}, p.77; this problem was also considered in
  \cite{St2}).   We will not pursue this here.

\section{Final remarks}\label{rela}

Our method should apply without major difficulties to some settings which had already been treated in the range $1<p<2$ 
by the use of the method in \cite{CD}:  graphs (\cite{R}) and vector bundles (\cite{TW}).  In the first
direction, some results already exist in the range $1<p<\infty$ in the group case, see  \cite{ND2}, \cite{I2}.
It should also be possible to cover the general setting of  \cite{Siko}.

Other directions are left open by the present work. Since we work in the framework of singular integrals theory (and, in
particular, use a
 volume growth assumption), the estimates we obtain do depend on the dimension, contrary
to the ones in
\cite{Ba0},
\cite{Ba}. To obtain dimension-free estimates, or to work in an infinite dimensional setting, is a
subject in itself, and was achieved so far only in rather specific situations: see  \cite{SteinR}, \cite{PAM},
\cite{P}, \cite{CMZ}, \cite{DV}, \cite{Pq}, and the references therein. 
Let us  point out that this is an advantage of the approach to Riesz transform boundedness
via Littlewood-Paley theory, as in  \cite{Ba0},
\cite{Ba}. See also \cite{CDfull}. We do not see for the time being how to cumulate the advantages
of both methods, which would yield a proof of the conjecture in \cite{CDfull}.

We do not touch either  the question of higher order Riesz transforms, whose boundedness would
require much more regularity; recall that already for Lie groups with polynomial growth, even the
$L^2$ boundedness of second order Riesz transforms only takes place on nilpotent groups and
their compact extensions (see \cite{A}, \cite{ERS}).

\bigskip

{\bf Acknowledgements:}   The authors would like to thank Marc Arnaudon,  Gilles Carron, Ewa Damek, Brian Davies, Li Hong-Quan,
Li Xiang-Dong, No\"el Lohou\'e, Alan McIntosh, Michel Rumin and Laurent Saloff-Coste for useful discussions,  references,
questions or remarks. Thanks are also due to S\"onke Blunck and  Ishiwata Satoshi for reading carefully the manuscript.

\bigskip

\end{document}